\documentclass[12pt]{article}
\usepackage{latexsym,amsmath,amssymb}

\newif\ifdviwin

\usepackage[latin1]{inputenc}
\usepackage[english]{babel}
\usepackage{indentfirst}
\usepackage[mathscr]{eucal}
\usepackage{amssymb,amsmath,amsfonts}
\usepackage{fancybox,fancyhdr}
\usepackage{graphicx}

\newif\ifdviwin

\dviwintrue

\def\cA{\mathcal{A}}

\def\cS{\mathcal{S}}
\def\cC{\mathcal{C}}

\def\cB{\mathcal{B}}

\def\cL{\mathcal{L}}
\def\cU{\mathcal{U}}
\def\cG{\mathcal{G}}
\def\cQ{\mathcal{Q}}
\def\cF{\mathcal{F}}

\let\8=\infty \let\0=\emptyset

\let\tilde=\widetilde

\let\landa=\lambda
\let\alfa=\alpha
\let\ro=\rho
\let\ex=\wedge

\let\parc=\partial

\def\landa{\lambda}

\def\lap{\Delta}
\def\flecha{\rightarrow}
\def\esiz{\langle}
\def\esde{\rangle}

\def\cte.{\mathop{\rm cte.}\nolimits}

\def\det{\mathop{\rm det}\nolimits}

\def\Re{\mathop{\rm Re }\nolimits}

\def\cosh{\mathop{\rm cosh }\nolimits}

\def\N{\mathbb{N}}
\def\L{\mathbb{L}}

\def\Z{\mathbb{Z}}
\def\Q{\mathbb{Q}}
\def\R{\mathbb{R}}

\def\C{\mathbb{C}}

\def\H{\mathbb{H}}
\def\S{\mathbb{S}}

\def\SL{\mathbf{SL}(2,\mathbb{C})}
\def\He{\rm{Herm}(2)}

 \newtheorem{defi}{Definition}
 \newtheorem{teo}[defi]{Theorem}
 
 \newtheorem{cor}[defi]{Corollary}
 
 \newtheorem{eje}[defi]{Example}
 \newtheorem{remark}[defi]{Remark}

 \newenvironment{proof}{\rm \trivlist \item[\hskip \labelsep{\it
      Proof}:]}{\par\nopagebreak \hfill $\Box$ \endtrivlist}

\numberwithin{equation}{section}

\begin{document}
\mbox{}\vspace{0.4cm}

\begin{center}
\rule{14cm}{1.5pt}\vspace{0.5cm}

\renewcommand{\thefootnote}{\,}
{\Large \bf The Cauchy problem for Liouville equation \\[0.3cm] and Bryant
surfaces}\\ \vspace{0.5cm} {\large José A. Gálvez$^a$ and Pablo Mira$^b$}\\
\vspace{0.3cm} \rule{14cm}{1.5pt}
\end{center}
  \vspace{1cm}
$\mbox{}^a$ Departamento de Geometría y Topología, Universidad de Granada,
E-18071 Granada, Spain. \\ e-mail: jagalvez@ugr.es \vspace{0.2cm}

\noindent $\mbox{}^b$ Departamento de Matemática Aplicada y Estadística,
Universidad Politécnica de Cartagena, E-30203 Cartagena, Murcia, Spain. \\
e-mail: pablo.mira@upct.es \vspace{0.3cm}

 \begin{abstract}
We give a construction that connects the Cauchy problem for Liouville
elliptic equation with a certain initial value problem for mean curvature
one surfaces in hyperbolic $3$-space $\H^3$, and solve both of them. We
construct the only mean curvature one surface in $\H^3$ that passes
through a given curve with given unit normal along it, and provide diverse
applications. In particular, topics like period problems, symmetries,
finite total curvature, planar geodesics, rigidity, etc. of surfaces are
treated.
 \end{abstract}

\section{Introduction}
The classical elliptic Liouville equation \cite{Lio}
 \begin{equation}\label{lio}
\lap (\log \phi ) = -2c \phi, \hspace{0.5cm} c\in \R
 \end{equation}
is an important research topic in partial differential equations, as
evidenced by the amount of works that has generated. It also has a clear
connection with differential geometry, since its solutions describe the
conformal factor that turns a flat metric into a pseudo-metric of constant
curvature $c$ on a surface. This interpretation shows that Liouville
equation admits a holomorphic resolution, namely, any solution to
\eqref{lio} on a simply connected domain
$\Omega\subseteq\C$ is of the form
 \begin{equation}\label{sollio}
\phi(s+it) = \frac{4|g'(z)|^2}{(1+c|g(z)|^2)^2},\hspace{1cm} z=s+it
 \end{equation}
where $g$ is an arbitrary meromorphic function on $\Omega$ (holomorphic
with $1>-c|g|^2$ if $c\leq 0$). In addition, the connection of Liouville
equation for
$c=1$ with minimal surface theory is well-known, and appears implicitly in
books as \cite{DHKW,Nit,Oss}.

A less typical interrelation of Liouville equation and surface theory
occurs in the study of surfaces with constant mean curvature one in the
standard $3$-dimensional hyperbolic space $\H^3$ (see for instance
\cite{Ten}). This was used by R. Bryant in his 1987's seminal paper
\cite{Bry}, in which he derived a holomorphic representation for these
surfaces, analogous in its spirit to the classical Weierstrass
representation of minimal surfaces in
$\R^3$. After Bryant's work, the above class of surfaces has become a
fashion research topic, and has received many important contributions
\cite{UY1,UY2,RUY1,RUY2,RUY3,Sma,CHR,HRR,Yu}. We shall use the term
\emph{Bryant surfaces} when referring to surfaces with mean curvature one
in $\H^3$.

Let us introduce the following \emph{Cauchy problem} for the class $\cB$
of Bryant surfaces:

 \begin{quote}
Let $\beta:I\flecha \H^3$ be a regular analytic curve, and $V:I\flecha
\S_1^3$ an analytic vector field along $\beta$ such that $\esiz
\beta,V\esde \equiv \esiz \beta', V\esde \equiv 0$. Find all Bryant
surfaces containing
$\beta$ and with unit normal in $\H^3$ along $\beta$ given by $V$.
 \end{quote}

The first objective of this paper is to give a back-and-forth construction
connecting the Cauchy problem for Liouville equation, and the Cauchy
problem for Bryant surfaces, and to solve both of them. The second
objective is to apply the solution of the Cauchy problem for
$\cB$ to study the geometry of Bryant surfaces.

We remark that the above formulated Cauchy problem for Bryant surfaces has
been inspired by the classical Björling problem for minimal surfaces in
$\R^3$, proposed by E.G. Björling in 1844 and solved by H.A.
Schwarz in 1890. Further details, as well as some research on this topic
may be consulted in \cite{DHKW,Nit,ACM,GaMi2,MiPa}.

We have organized this paper as follows. In Section 2 we study the Cauchy
problem for Liouville equation, by using both analytic and geometric
methods. For the analytic part, and inspired by \cite{Lio}, we view
\eqref{lio} as a complex differential equation, and solve the Cauchy
problem for it by means of another classical tool, the \emph{Schwarzian}
derivative. Geometrically, we show that solving the Cauchy problem for
\eqref{lio} is equivalent to the problem of integrating the Frenet
equations for curves in the standard $2$-dimensional Riemannian space
model $\cQ(c)$ of constant curvature $c\in \R$. This interpretation
provides an explicit resolution of the Cauchy problem for \eqref{lio} and
certain adequate initial data, which will be of a special interest when
regarded in the context of Bryant surfaces. In addition, we also solve
explicitly the Cauchy problem for the degenerate version of Liouville
equation, i.e. $\lap (\log \phi)=0$. This corresponds to integrate the
Frenet equations in $\R^2$, which is a typical exercise for undergraduate
students.

In Section 3 we show that, given a pair of Björling data $\beta,V$, there
is a unique solution to the Cauchy problem for Bryant surfaces with
initial data $\beta,V$, and we construct such Bryant surface in terms of
the solution of the Cauchy problem for Liouville equation with $c=1$ given
in Section 2. This construction is essentially self-contained, as it does
not use the Bryant representation in \cite{Bry}. The relation of Liouville
equation with Bryant surfaces comes from the following fact: if
$\psi:\Sigma\flecha \H^3$ is a Bryant surface with metric $ds^2$ and curvature $K$ ($\leq 0$), then
$-K ds^2$ is a pseudo-metric of constant curvature one, i.e. if $z$
is a conformal parameter of $\Sigma$, then $-Kds^2 = \phi |dz|^2$ and
$\phi$ verifies Liouville equation for $c=1$.  We also present in Section 3,
this time using the Bryant representation, a simplified construction of
the only solution to the Cauchy problem for
$\cB$ by means of two important equations of the theory, due to
Umehara-Yamada \cite{UY1} and Small \cite{Sma}, respectively.

In Section 4 we regard the solution of the Cauchy problem for $\cB$ as a
meromorphic representation for Bryant surfaces, and explore its
applications. Before describing the results obtained there, we note the
following facts on Bryant surfaces.

\begin{enumerate}
\item
For a Bryant surface, the developing map $g$ of the pseudo-metric $-K
ds^2$ is generally not single-valued on the surface. As a consequence, the
total curvature of a complete Bryant surface does not admit a
quantization, and can assume any negative value. These limitations do not
appear in minimal surface theory.
 \item
To construct non simply-connected minimal surfaces in $\R^3$, one only
needs to ensure that a certain holomorphic differential on a Riemann
surface has no real periods (what, obviously, can be difficult). In
contrast, the period problem for Bryant surfaces is a more complicated
question, which difficulties the construction of Bryant surfaces with
non-trivial topology.
 \item
While in $\R^3$ the coordinates of a minimal surface are obtained in terms
of the Weierstrass data $(g,\omega)$ after computing an integral, in
Bryant surface theory the explicit coordinates are only obtained from
$(g,\omega)$ after solving a second order differential equation. As this equation
cannot be explicitly integrated except for some special cases, it is not
known how to construct Bryant surfaces in explicit coordinates with
prescribed geometrical properties. This explains why some explicit
formulas in minimal surface theory, such as the associate family transform
and the López-Ros transform (a conformal method to deform a minimal
surface so that one of its coordinates remains invariant) have not been
extended to Bryant surfaces.
\end{enumerate}
In Section 4 we treat these questions by means of the solution to the
Cauchy problem for Bryant surfaces. Regarding the first point, we give a
criterion to determine when the developing map $g$ is single-valued on a
Bryant surface, in terms of the following problem: when is a curve in
$\S^2$ with periodic geodesic curvature closed? Regarding period problems, we
give a symmetry principle for Bryant surfaces and show as a corollary that
the meromorphic representation yielded by our solution to the Cauchy
problem for $\cB$ provides a period-problem free representation of Bryant
cylinders. As applications, we give a general description of the complete
Bryant surfaces with genus zero, two ends, and that have finite total
curvature of finite dual total curvature. Important results in this
direction have been obtained in \cite{UY1,RUY2,RUY3}. Concerning the third
point, we describe a method with constructs in explicit coordinates the
only two Bryant surfaces that contain a given plane curve as a planar
geodesic. We also show in Section 4 how the solution to the Cauchy problem
for $\cB$ provides a simple classification of Bryant surfaces which are
invariant under $1$-parameter isometry groups of $\H^3$. In particular, we
prove without considering differential equations (see \cite{Nit} in
contrast for the case of minimal surfaces in $\R^3$) that any helicoidal
Bryant surface is associate to a rotational Bryant surface. This line of
inquiry has been motivated by \cite{EaTo}, in where it is shown that the
converse of this result is not true, i.e. not all Bryant surfaces
associated to a rotational one are helicoidal.

Finally, in Section 5 we use our analysis on the Cauchy problem for Bryant
surfaces to provide a more explicit geometric solution of the Cauchy
problem for Liouville equation (for $c=1$), and to investigate the
integration of the Frenet equations for curves in $\S^2$.

This work is part of the PhD Thesis of the second author, which was
defended at the University of Murcia in 2003.

\section{Liouville equation}

Let $\cU\subseteq\C$ be a planar domain, and consider the usual Wirtinger
operators $\parc_z=(\parc_s -i\parc_t )/2$, and
$\parc_{\bar{z}}=(\parc_s+i\parc_t)/2$, being
$z=s+it$. Then Liouville equation is written as
 \begin{equation}\label{licom}
 (\log \phi)_{z\bar{z}} = -(c/2) \phi,
 \end{equation}
and once in this form can be considered as a complex differential
equation. If $ds^2=\landa |dz|^2$ is a Riemannian metric on $\cU$, then
$\landa$ satisfies Liouville
equation \eqref{licom} if and only if $ds^2$ has constant curvature of
value $c$. Our first aim in this Section is to provide a purely analytic
resolution of the Cauchy problem for Liouville equation:
 \begin{equation}\label{caulio}
\left\{\def\arraystretch{1.2} \begin{array}{lll}\displaystyle \left(\log
\phi\right)_{z\bar{z}}&=&\displaystyle -\frac{c}{2}\phi, \\ \phi(s,0)&=&
a(s),\\ \phi_z(s,0)&=& b(s).
\end{array}\right.
 \end{equation}
Here, in order for the solution to be real, we ask $a(s)$ to be a
non-negative real analytic function, and $b(s)$ to be a real analytic
complex function such that $2 \Re b(s) = a'(s)$. Actually, we will always
assume for clarity that $a(s)$ is positive. This makes no restriction,
since this condition can be suppressed \emph{a posteriori} by analytic
continuation.

The solution we are going to expose relies on finding a complex-valued
function with prescribed Schwarzian derivative. Recall that the
\emph{Schwarzian derivative} of a meromorphic function
$f$ with respect to a complex parameter $z$ is $$\{f,z\}=
\left(\frac{f''}{f'}\right)' - \frac{1}{2}\left(\frac{f''}{f'}\right)^2,
\hspace{1cm} \left( '=\frac{d}{dz}\right).$$ The same definition applies
when we consider $f$ to be a real function, and instead of the complex
parameter $z$ we write a real parameter $s$.

\begin{teo}\label{calite}
Let $c\neq 0$, and $T(s):I\flecha \C$ be an arbitrary solution of the
differential equation $\{T,s\}= \Upsilon(s)$, where
 \begin{equation}\label{ucabr}
2\Upsilon(s)=
2\left(\frac{b(s)}{a(s)}\right)'-\left(\frac{b(s)}{a(s)}\right)^2 +c \,
a(s),
  \end{equation}
and define
 \begin{equation}\label{rcabr}
\overline{R(s)}=\frac{1}{c} \,\frac{T''(s) -(b(s)/a(s)) T'(s)}{2 T'(s)^2 -
T(s) T''(s) +(b(s)/a(s)) T(s) T'(s)}.
 \end{equation}
If $T(z),R(z)$ are meromorphic extensions of
$T(s),R(s)$ on an open subset $D\subset\C$ containing $I$,
then
 \begin{equation}\label{ficabr}
\phi(s,t)= \frac{4 T_z\overline{R_z}}{(1+c T\bar{R})^2}
 \end{equation}
is the only solution to the Cauchy problem \eqref{caulio}.
\end{teo}
\begin{proof}
We begin by noting that if $T,R$ are arbitrary meromorphic functions and
$z=s+it$, then the map \eqref{ficabr} satisfies Liouville equation
\eqref{licom}.

Now, let $U(s),V(s)$ be two functions from $I$ into $\C$ such that $U^2=
R'$ and $V^2=T'$, and define

$$\def\arraystretch{2.5}\begin{array}{ll}
L(s)=\displaystyle\frac{1}{U (s)}, & \hspace{0.4cm}
H(s)=\displaystyle\frac{\sqrt{c}\, R(s)}{U(s)} ,
\\ M(s)=\displaystyle\frac{1}{V(s)}, & \hspace{0.4cm}
N(s)=\displaystyle\frac{\sqrt{c}\, T(s)}{V(s)}
\end{array} $$ where $\sqrt{c}\in \C$ is a fixed square root of $c$.
From here and \eqref{ficabr}, denoting $\varepsilon ={\rm sign} (c)$,
we see that
 \begin{equation}\label{film}
 \phi(s,0)=\frac{4}{\left(\bar{L}(s) M(s) +\varepsilon \bar{H}(s)
 N(s)\right)^2}.
 \end{equation}
In addition, differentiation of \eqref{rcabr} yields
$$ \bar{R}' = \frac{1}{c} \frac{{T'}^3 \big(2\{T,s\} +(b/a)^2
-2(b/a)'\big)}{\left(2 {T'}^2 -T T'' + (b/a) T T'\right)^2},$$ where we
have suppressed the parameter $s$. Making use of \eqref{ucabr}, we obtain
from this expression that
 \begin{equation}\label{rcapri}
 \delta {\bar U} = \frac{\sqrt{a} T' V}{2 {T'}^2 -T T'' + (b/a)
T T'},
 \end{equation}
where here $\delta =\pm 1$. Now, from \eqref{rcapri} it is direct to check
the relation
 \begin{equation}\label{lebar}
 \bar{L}=\frac{2\delta }{\sqrt{c} \,\sqrt{a}} \left(N'+\frac{b}{2a} N\right).
 \end{equation}
In the same manner, we may express $H$ in terms of $M$ as
  \begin{equation}\label{hbal}
 \bar{H}=\frac{-2\varepsilon \delta}{\sqrt{c} \, \sqrt{a}} \left(M'+\frac{b}{2a} M
 \right).
 \end{equation}
Besides, since $N'M-NM'=\sqrt{c}$, equations \eqref{lebar} and
\eqref{hbal} let us obtain $(\bar{L}M +\varepsilon \bar{H} N)^{-2} = a/4$.
Therefore, \eqref{film} indicates that $\phi(s,0)=a(s)$ on $I$.

To check that the remaining initial condition is fulfilled, we first note
that, from \eqref{ficabr},
 $$\phi_z (s,t)=\frac{4}{(1+c T\bar{R})^3} \left\{T_{zz} \bar{R_z} (1+c T\bar{R})
 -2 c T_z^2 \overline{R R_z}
\right\}.$$ Then, since $T,R$ are meromorphic, we get $$ \phi_z
(s,0)=\frac{4}{(1+c T\bar{R})^3} \left\{T'' \bar{R'} (1+c T\bar{R}) -2 c
{T'}^2 \overline{R R'} \right\}.$$ Observing now that
$\bar{L}M +\varepsilon \bar{H}N = (1+ c T\bar{R})/(V \bar{U})$,
a direct computation ensures that $$\phi_z (s,0)=\frac{4}{(\bar{L} M
+\varepsilon \bar{H}N)^3} \left\{\frac{T''(1+c T\bar{R}) -2c {T'}^2
\bar{R}}{\bar{U} V T'}\right\} =\frac{-8\left( \bar{L}M' +\varepsilon
\bar{H} N'\right)}{(\bar{L} M +\varepsilon \bar{H}N)^3}.$$ From this
expresion, and since
$\bar{L}M +\varepsilon \bar{H} N =2\delta / \sqrt{a}$, equations \eqref{lebar} and
\eqref{hbal} let us conclude that $\phi_z(s,0) = b(s)$. Hence,
$\phi(s,t)$ is the solution to the Cauchy problem \eqref{caulio}.
 \end{proof}
\begin{remark}
Some of the functions we have introduced are not well defined at certain
points. However, the proof keeps working if we skip those points and use
analytic continuation afterwards. We have decided not to make this detail
explicit at every stage of the proof in order to gain clarity.
\end{remark}

Next, we study geometrically the Cauchy problem \eqref{caulio}, this time
regarded in a real form, i.e.
 \begin{equation}\label{caulior}
\left\{\def\arraystretch{1.2} \begin{array}{lll}\displaystyle \lap (\log
\phi)&=&\displaystyle -2c \phi, \\ \phi(s,0)&=& a(s),\\ \phi_t(s,0)&=&
d(s).
\end{array}\right.
 \end{equation}
Here $a(s),d(s)$ are real analytic functions, and $a(s)$ is positive.

To do so, we denote by $\cQ (c)$ the standard $2$-dimensional Riemannian
space form of constant curvature $c$, that is,
$\cQ(0)=\R^2$ and
 \begin{equation*}
\cQ (c)= \left\{ \def\arraystretch{2.5}
\begin{array}{ll}\displaystyle \S^2(c)= \left\{(x_0,x_1,x_2):
x_0^2+x_1^2+x_2^2= \frac{1}{\sqrt{c}} \right\} & \text{if } c>0, \\
\displaystyle \H^2(c)= \left\{(x_0,x_1,x_2): -x_0^2+x_1^2+x_2^2=
\frac{-1}{\sqrt{-c}}, \ x_0>0 \right\} & \text{if } c<0.
\end{array}
\right.
 \end{equation*}
We shall also consider the \emph{stereographic projection} $\pi$ of $\cQ
(c)$ into $\C\cup \{\8\} $, defined by $$ \pi (x_0,x_1,x_2)= \frac{x_1+i
x_2}{1-c x_0} .$$ With this, we have the following geometric description
of the solution to the Cauchy problem \eqref{caulior}.

\begin{teo}\label{curvas}
The only solution to the Cauchy problem for Liouville equation
\eqref{caulior} is $$\phi(s,t)= \frac{4 |g'(z)|^2}{(1+c |g(z)|^2)^2},
\hspace{1cm} z=s+it,$$ where $g(z)$ is the meromorphic extension of
$g(s)=\pi (\alfa(s))$, being $\pi:\cQ(c) \flecha \C\cup \{\8\}$ the stereographic
projection, and $\alfa(s)$ the only curve in $\cQ(c)$ with arclength
parameter and geodesic curvature given respectively by
 \begin{equation}\label{paracurv}
 u(s)= \int^s \sqrt{a}(r) dr, \hspace{1cm} \text{and} \hspace{1cm}
 \kappa(s)= \frac{-d(s)}{2 a(s)^{3/2}}.
 \end{equation}

\end{teo}
\begin{proof}
Let $\phi(s,t)$ denote the solution to \eqref{caulior}, defined on a
simply connected complex domain $\Omega\subseteq \C$ containing $I$. Then,
from the Frobenius theorem we can assure the existence of a conformal map
$F:\Omega\subseteq\C\flecha \cQ(c)$ such that $\esiz dF,dF\esde = \phi
|dz|^2$. If we denote $\alfa(s)= F(s,0)$, this relation indicates that
$a(s)=\esiz \alfa'(s),\alfa'(s)\esde$. In particular, $\alfa(s)$ is a
regular curve in $\cQ(c)$, and its arclength parameter is the one
specified in \eqref{paracurv}.

Let now $J$ denote the complex structure on $\cQ(c)$, and $\kappa(s)$ the
geodesic curvature of $\alfa(s)$. Then
$$
\def\arraystretch{2}\begin{array}{lll}
 d(s) & = & \phi_t (s,0) = \displaystyle\frac{\parc}{\parc t} \esiz F_s,F_s \esde (s,0) = 2
 \left\esiz \displaystyle\frac{\parc}{\parc s} F_t, F_s\right\esde(s,0) \\ & = & 2 \left\esiz
 \displaystyle\frac{d}{ds}J \alfa'(s),\alfa'(s)\right\esde = -2 |\alfa'(s)|^3 \kappa (s).
\end{array} $$
Once here, the proof concludes by a standard application of the identity
principle for holomorphic functions.
\end{proof}

This Theorem shows in particular that the choice $d(s)=0$ corresponds to
geodesics in $\cQ(c)$. Thus, in the specific case of $c=1$ (the one we
shall use in the geometric part of this paper), we find the following.
 \begin{cor}\label{corgp}
Let $a(s):I\flecha \R^+$ be real analytic, and $z=s+it$. The only solution
to the Cauchy problem
\begin{equation}\nonumber\left\{\def\arraystretch{1.2} \begin{array}{lll}\displaystyle  \lap \log
\phi &=&\displaystyle -2 \phi, \\ \phi(s,0)&=& a(s),\\ \phi_t(s,0)&=& 0
\end{array}\right.
\end{equation}
is
 \begin{equation}\nonumber
\phi(s,t)=\frac{4|g_z|^2}{(1+|g|^2)^2}, \hspace{0.5cm} \text{being}
\hspace{0.5cm} g(z)=\exp\left( i\int_{s_0}^z \sqrt{a}(\zeta )
d\zeta\right).
 \end{equation}
Here $\sqrt{a}(z)$ is a holomorphic extension of $\sqrt{a(s)}$.
 \end{cor}

Rather than with Liouville equation, we shall treat the Cauchy problem for
Bryant surfaces with the \emph{modified Liouville equation}
 \begin{equation}\label{limod}
 4(\log \rho)_{z\bar{z}} = -\rho^2 |f(z)|^2,
 \end{equation}
where $f$ is a meromorphic function. The relation of this equation with
the Liouville one \eqref{licom} is very tight, as if
$\rho:D\subseteq\C\flecha \R$ is smooth, then $\ro$
satisfies \eqref{limod} if and only if $\phi=\ro^2 |f|^2$ satisfies
\eqref{licom} for $c=1$.

This comments are enough to solve the Cauchy problem associated to
\eqref{limod},
\begin{equation}\label{caulimo}\left\{\def\arraystretch{1.2} \begin{array}{lll}\displaystyle  4\left(\log
\ro\right)_{z\bar{z}}&=&-\ro^2 |f(z)|^2,
\\ \ro(s,0)&=& v(s),\\ \ro_z(s,0)&=& w
(s).
\end{array}\right.
\end{equation}
Here, we assume that $f$ does not have poles on $\R$.
 \begin{cor}\label{sollimo}
The solution to the Cauchy problem for the modified Liouville equation
\eqref{caulimo} is
$\ro=\sqrt{\phi}/ |f|$, being $\phi$ the solution to the Cauchy problem
for Liouville equation for $c=1$ and the initial data
 \begin{equation}\label{insolli}
 \def\arraystretch{1.2}\begin{array}{lll}
 a(s) & = & v(s)^2 |f(s)|^2, \\
 b(s) & = & 2 v(s) w(s) |f(s)|^2 + v(s)^2 f'(s)\overline{f(s)}
 \end{array}
 \end{equation}
constructed in Theorem \ref{calite}.
\end{cor}

The degenerate case of Liouville equation, i.e. the case in which
$c=0$, is not covered by  Theorem \ref{calite}. However, in this particular case the
Cauchy problem
\begin{equation}\label{intro4}
\left\{\def\arraystretch{1.2} \begin{array}{lll}\displaystyle \lap
\left(\log \phi\right)&=& 0, \\ \phi(s,0)&=& a(s),\\ \phi_t(s,0)&=& d(s),
\end{array}\right.
\end{equation}
admits an explicit holomorphic resolution.

 \begin{teo}\label{solide}
Let $a(s),d(s):I\flecha \R$ be real analytic functions with
$a(s)$ positive. The only solution $\phi$ to the Cauchy problem \eqref{intro4}
is constructed as follows: choose $s_0\in I$ arbitrary, let
$\theta(s):I\flecha \R$ be
 \begin{equation}\label{detita}
\theta(s)=-\frac{1}{2} \int_{s_0}^s \frac{b(r)}{a(r)} dr,
 \end{equation}
and take holomorphic extensions $\sqrt{a}(z)$,
$\theta(z)$ of
$\sqrt{a(s)}$, $\theta(s)$. Then $\phi(s,t)=|f(z)|^2$, where $z=s+it$ and
$f$ is the holomorphic function
 \begin{equation}\nonumber
 f(z)=\sqrt{a}(z) e^{i\theta(z)}.
 \end{equation}
 \end{teo}
\begin{proof}
We give a constructive proof of this result. First, observe that any
solution to $\lap (\log \phi )=0$ on a simply connected domain $U$ is of
the form $\phi=|f|^2$ for some holomorphic function $f$ on $U$. If
$\phi$ solves the Cauchy problem \eqref{intro4}, and $\theta(s)$ denotes the argument
of $f$ along $I$, then
 \begin{equation}\label{efceta}
 \def\arraystretch{1.6}\begin{array}{lll}
f(s,0) & = & \sqrt{a(s)} e^{i\theta(s)}, \\ f_z(s,0) & = & \left\{
\left(\sqrt{a(s)}\right)' + i \theta'(s) \sqrt{a(s)}\right\}
e^{i\theta(s)}.
 \end{array}
 \end{equation}
In addition, we have
$$\frac{1}{2}\left(a'(s)-ib(s)\right) =\phi_z(s,0)= f_z(s,0) \overline{f(s,0)},$$
from where it is obtained using \eqref{efceta} that
$$\frac{1}{2}\left(a'(s)-ib(s)\right) = \frac{1}{2} a'(s)
+i\theta'(s) a(s),$$ i.e. $\theta(s)$ is given by \eqref{detita}. Hence,
$f(s,0)= \sqrt{a(s)}e^{i\theta(s)}$ for all $s\in I$, and the proof concludes by analytic continuation.
\end{proof}

According to Theorem \ref{curvas}, solving the Cauchy problem
\eqref{intro4} is equivalent to the problem of integrating the Frenet
equations for curves in $\R^2$. However, the solution to this second
problem is well known, and appears in most textbooks on classical
differential geometry. Specifically, the result states that the only (up
to congruences) curve in $\R^2$ parametrized by arclength and with
prescribed curvature $k(s)$ is $$\alfa(s)= \left( \int^s \cos \theta(u)
du, \int^s \sin \theta(u) du\right), \hspace{1cm} \theta(s)=\int^s k(u)
du.$$ It is then possible to provide an alternative geometric proof of
Theorem \ref{solide} by means of this equation and Theorem \ref{curvas}.

We also remark that the equation $\lap (\log \phi) =0$ plays an important
role in the study of flat surfaces in the hyperbolic
$3$-space \cite{GMM1,GMM2,GaMi3,KUY}.

\section{The Cauchy problem for Bryant surfaces}

We begin by recalling some basic facts of the Hermitian model for the
hyperbolic $3$-space. So, let $\L^4$ denote the
$4$-dimensional Lorentz-Minkowski space, that is, the real vector space
$\R^4$ endowed with the Lorentzian metric
$$\esiz,\esde=-dx_0^2+dx_1^2+dx_2^2+dx_3^2,$$ in canonical coordinates.
The cross product of $u_1,u_2,u_3\in \L^4$ is defined by the identity
$\esiz u_1\times u_2\times u_3 , w\esde = {\rm det} (u_1,u_2,u_3,w)$.
We shall identify $\L^4$ with the space of $2$ by
$2$ Hermitian matrices in the usual way,
$$(x_0,x_1,x_2,x_3)\in\L^4\longleftrightarrow
\left(\begin{array}{cc} x_0+x_3 & x_1+i x_2\\ x_1-i x_2& x_0-x_3\\
\end{array}\right)\in {\rm Herm }(2).$$ Under this identification one gets
$\esiz m,m\esde=-\det (m)$ for all $m\in {\rm Herm}(2)$. The complex Lie group $\SL$
acts on $\L^4$ by $\Phi \cdot m =\Phi m \Phi^*$, being $\Phi\in \SL$,
$\Phi^*=\bar{\Phi}^{t}$, and $m\in \He$. Consequently, $\SL$ preserves the metric
and the orientations. We shall realize the hyperbolic $3$-space of
negative curvature $-1$ in its Minkowski model, that is, $\H^3=\{x\in
\L^4: \esiz x,x\esde=-1, x_0>0\}.$ Observe that, for every $q\in \H^3$ and
every $u_1,u_2 \in T_q \H^3$, we can define the \emph{exterior product} as
$u_1 \wedge u_2 = p\times u_1 \times u_2$.

The above identification makes $\H^3$ become
$$\H^3= \left\{ \Phi \Phi^* : \Phi \in \SL\right\}.$$
In the same way, the \emph{positive light cone} $\N^3=\{x\in \L^4 : \esiz
x,x\esde=0, x_0>0 \}$ is seen as the space of positive semi-definite
matrices in $\He$ with determinant $0$, and can be described as
$$\N^3= \left\{ w w^* : w=(w_1,w_2)\in \C^2\right\},$$ where $w\in \C^2$
is uniquely defined up to multiplication by an unimodular complex number.
The quotient $\N^3/\R^+$ inherits a natural conformal structure and it can
be regarded as the ideal boundary $\S_{\8}^2 $ of the hyperbolic $3$-space
$\H^3$ in $\L^4$. The map $ww^*\flecha [(w_1,w_2)]$ becomes the quotient
map of $\N^3$ onto $\S_{\8}^2$ and identifies $\S_{\8}^2$ with
$\C\mathbf{P}^1$. By stereographic projection, we may also regard this
quotient map as $w w^* \flecha w_2 / w_1$ from $\N^3$ into $\C\cup \{\8\}
$.

Let $\psi:\Sigma\flecha\H^3$ be an immersed surface in $\H^3$, with unit
normal $\eta:\Sigma\flecha \S_1^3$ in $\H^3$. Here, $\S_1^3=\{x\in \L^4:
\esiz x,x\esde =1\}$ is the \emph{de Sitter} $3$-space. We shall regard
$\Sigma$ as a Riemann surface with the conformal structure induced by
the isometric immersion
$\psi$. Then, we can consider the hyperbolic Gauss map of the surface, $G:\Sigma\flecha \C\cup \{\8\}$,
defined by $G=[\psi+\eta]$. In other words, if we denote $N= \psi +\eta$,
then $G$ may be regarded as the map
 \begin{equation}\label{gamacu}
G= \frac{N_1-iN_2}{N_0+N_3}: \Sigma\flecha \C\cup \{\8\}, \hspace{0.7cm}
N=(N_0,N_1,N_2,N_3).
 \end{equation}
Bryant proved in \cite{Bry} the fundamental fact that $G$ is meromorphic
if and only if the surface has mean curvature one, $H\equiv 1$. As we have
already mentioned, we will call such surfaces \emph{Bryant surfaces}. We
also remark that the \emph{Hopf differential}, defined by $Q= \esiz
\psi_{z z}, \eta\esde dz^2$, is a globally defined \emph{holomorphic}
$2$-form on
$\Sigma$ whenever the immersion $\psi$ is a Bryant surface.

To solve the Cauchy problem for Bryant surfaces specified in Section one,
we shall call any pair $\beta,V$ in the conditions of that problem a pair
of \emph{Björling data}. Then we have:

 \begin{teo}\label{solbjbr} Given Björling data $\beta(s),V(s)$, there is a unique solution
to the Cauchy problem for Bryant surfaces with initial data
$\beta(s),V(s)$. This solution
$\psi:D\subseteq\C\flecha \H^3$ can be constructed in a neighbourhood of $\beta$ as follows:
let $\beta(z),V(z)$ be holomorphic extensions of
$\beta(s),V(s)$, and define $\nu(z)=\beta(z)+V(z)$ and
 \begin{equation}\label{himabj}
 G(z)=\frac{\nu_1(z)-i\nu_2(z)}{\nu_0(z)+\nu_3(z)}.
 \end{equation}
Let $\ro:D\subseteq\C\flecha [0,+\8)$ be the only solution to the Cauchy
problem
\begin{equation}\label{caulb}\left\{\def\arraystretch{1.2} \begin{array}{lll}\displaystyle  4\left(\log
\ro\right)_{z\bar{z}}&=&-\ro^2 |G_z|^2,
\\ \ro(s,0)&=& \nu_0(s)+\nu_3(s),\\
\ro_z(s,0)&=& \frac{1}{2}\left\{\nu_0'(s)+\nu_3'(s)
+i[\big(\nu(s)\wedge\nu'(s)\big)_0 +
\big(\nu(s)\wedge\nu'(s)\big)_3]\right\},
\end{array}\right.
\end{equation}
constructed via Corollary \ref{sollimo}. Then
$\psi=F\Omega F^*:D\subseteq\C\flecha \H^3$, where
 \begin{equation*}
 F= \left(\def\arraystretch{1.3}\begin{array}{cc} 1 & 0 \\
G & 1\\
\end{array}\right), \hspace{1cm}
\Omega= \left(\def\arraystretch{2}\begin{array}{cc} \ro
+\displaystyle\frac{2\ro_{z\bar{z}}}{\ro^2 |G_z|^2} & \displaystyle
\frac{2\ro_z}{\ro^2 G_z}
\\ \displaystyle
\frac{2\ro_{\bar{z}}}{\ro^2 \overline{G_z}} & \displaystyle\frac{2}{\ro}
\end{array}\right).
 \end{equation*}
Moreover, the hyperbolic Gauss map
$G:D\subseteq\C\flecha \C\cup \{\8\}$ of $\psi$ is given by
\eqref{himabj}, and its Hopf differential is
 \begin{equation}\label{jof}
 Q=-\frac{1}{2} \big\esiz \beta'(z)+V'(z),\beta'(z) -iV(z)\wedge
 \beta'(z)\big\esde dz^2.
 \end{equation}
 \end{teo}
\begin{proof}
First, we deal with the uniqueness part, and determine the form of the
solution.

Let $\beta(s),V(s)$ be Björling data defined on a real interval
$I$, and consider a Bryant surface that solves the Cauchy problem for these data.
Then, it is possible to parametrize conformally this surface in a
neighbourhood of $\beta$ as
$\psi:D\subseteq\C\flecha \H^3$, so that
 \begin{enumerate}
 \item
$D\subseteq\C$ contains an interval $J\subseteq I$.
 \item
$\psi(s,0)=\beta(s)$ for all $s\in J$, being $z=s+it$.
 \item
The unit normal to $\psi$ along $\beta(s)$ is $V(s)$.
 \end{enumerate}

Let $\eta:D\subseteq\C\flecha \S_1^3$ denote the unit normal to
$\psi$ in $\H^3$, and recall the hyperbolic Gauss map \eqref{gamacu}, defined in
$D$, where $N=\psi +\eta$. Consider also the curve $\nu(s):I\flecha \N^3$ given by
$\nu(s)=\beta(s) +V(s)$. Then
 \begin{equation}\label{genues}
 G(s)= \frac{\nu_1(s)-i\nu_2(s)}{\nu_0(s)+\nu_3(s)}.
 \end{equation}
Observe that, composing with a Möbius transformation in $\S_{\8}^2$ if
necessary, we may assume that
$\nu_0(s)+\nu_3(s)\neq 0$ for all $s\in I$.

Let $U\subseteq D$ be an open subset containing $J$ over which
$\beta(s),V(s)$ have holomorphic extensions $\beta(z),V(z)$. As $G$ is
meromorphic, we find from \eqref{genues} that
$G:U\subseteq\C\flecha \C_{\8}$ is given by \eqref{himabj}, and extends
meromorphically to all $D$.

In addition,
 \begin{equation}\label{bjh1}
 \psi +\eta =\ro \left(\def\arraystretch{1.1}\begin{array}{cc} 1 & \bar{G} \\
G & G \bar{G}\\
\end{array}\right),
 \end{equation}
being $\ro:D\flecha [0,+\8)$ the map $\ro=N_0+N_3$. But now, the fact that
$\psi$ is a Bryant surface provides \cite{Bry,GMM2} that $\esiz
d(\psi+\eta),d(\psi+\eta)\esde$ is a conformal pseudo-metric of curvature
one over $D$. Since \eqref{bjh1} yields
 \begin{equation}\label{bjh2}
\esiz (\psi+\eta)_z,(\psi+\eta)_{\bar{z}}\esde = \frac{1}{2}\ro^2 |G_z|^2,
 \end{equation}
we obtain that $\phi=\ro^2|G_z|^2$ satisfies Liouville equation
\eqref{licom} for $c=1$. In this way, $\ro:D\subseteq\C\flecha [0,+\8)$
satisfies the modified Liouville equation \eqref{limod}.

Observe next that, from $\ro=N_0+N_3$,
 \begin{equation}\label{bjh3}
\ro(s,0)=\nu_0(s)+\nu_3(s).
 \end{equation}
Besides, if $z=s+it$, it holds $\psi_t =\eta\wedge\psi_s$. As we saw that
$\psi +\eta$ is conformal, $(\psi+\eta)_s$ and
$(\psi+\eta)_t$ are orthogonal and with the same length. Moreover, the
basis $\{(\psi+\eta)_s,(\psi+\eta)_t\}$ is negatively oriented on the
(oriented) tangent plane of the immersion \cite{Bry}. Consequently,
$(\psi+\eta)_t = -\eta\wedge (\psi+\eta)_s$ on $D$. As
 \begin{equation}\label{bjh4}
 \psi_t(s,0)=V(s)\wedge \beta'(s),
 \end{equation}
the above relation tells that
 \begin{equation}\label{bjh5}
 -\eta_t(s,0) = V(s)\wedge \big(2\beta'(s) + V'(s)\big).
 \end{equation}
Equations \eqref{bjh4}, \eqref{bjh5} let us recover in terms of the
Björling data the Hopf differential
$Q= q(z) dz^2$, being $q(z)=\esiz \psi_{zz},\eta\esde$. Specifically, we obtain
\begin{equation}\label{cuese}
q(s)=-\frac{1}{2} \esiz \beta'+V', \beta' -iV\wedge \beta'\esde.
\end{equation}
Now, since $q(z)$ is holomorphic,
 \begin{equation}\label{bjh6}
q(z)= -\frac{1}{2} \big\esiz \beta'(z)+V'(z),\beta'(z) -iV(z)\wedge
 \beta'(z)\big\esde
 \end{equation}
on $D$. In addition, $(\psi+\eta)_t = -\eta\wedge (\psi+\eta)_s$ also
yields

\begin{equation}\nonumber
(\psi+\eta)_z(s,0)= \frac{1}{2}\big\{ \nu'(s) +i V(s)\ex \nu'(s) \big\}.
\end{equation}
Particularly, we conclude that
 \begin{equation}\label{bjh8}
\ro_z(s,0)= \frac{1}{2} \left\{\nu_0'(s)+\nu_3'(s)
+i\left[\big(V(s)\wedge\nu'(s)\big)_0 +
\big(V(s)\wedge\nu'(s)\big)_3\right]\right\}.
 \end{equation}
Hence, $\ro$ is the solution to the Cauchy problem for the modified
Liouville equation with initial data given by \eqref{bjh3} and
\eqref{bjh8}. Furthermore, this problem may be solved by means of
Corollary \ref{sollimo}, and $\ro$ is recovered in terms of the Björling
data.

Let us define next the meromorphic curve
$F:D\subseteq\C\flecha \SL$ given by
 \begin{equation}\label{bjh88}
 F= \left(\def\arraystretch{1.1}\begin{array}{cc} 1 & 0 \\
G & 1\\
\end{array}\right).\end{equation}
Then \eqref{bjh1} writes down as
 \begin{equation}\label{bjh9}
\psi +\eta = F\left(\def\arraystretch{1.1}\begin{array}{cc} \ro & 0
\\ 0 & 0\\
\end{array}\right) F^*.
 \end{equation}
Differentiation of this expression yields
  \begin{equation}\label{bjh10}
(\psi +\eta)_{\bar{z}} = F\left(\def\arraystretch{1.2}\begin{array}{cc}
\ro_{\bar{z}} & \ro \overline{G_z} \\ 0 & 0\\
\end{array}\right) F^*,
 \end{equation}
and

  \begin{equation}\label{bjh105}
(\psi +\eta)_{z\bar{z}} =F\left(\def\arraystretch{1.1}\begin{array}{cc}
\ro_{z\bar{z}} & \ro_z \overline{G_z} \\ \ro_{\bar{z}} G_z & \ro |G_z|^2\\
\end{array}\right)F^*.
 \end{equation}

On the other hand, let $\landa$ denote the conformal factor of the metric
of $\psi$, i.e. the positive smooth function such that $\esiz
d\psi,d\psi\esde = \landa |dz|^2$. Since the relation
 \begin{equation}\label{recgh}
 (\psi +\eta)_z = \frac{-2 q}{\landa}\psi_{\bar{z}}
 \end{equation}
holds on any Bryant surface, it is obtained that
$$(\psi+\eta)_{z\bar{z}} = -\frac{2\bar{q}}{\landa} \psi_{zz} +
\frac{2\bar{q} \landa_z}{\landa^2} \psi_z.$$ Besides, it is easy to check
the general relation
$$\psi_{zz}=\frac{\landa_z}{\landa} \psi_z +q \eta.$$
Putting together these two expressions we are left with
 \begin{equation}\label{bjh11}
\psi=\psi +\eta +\frac{\landa}{2|q|^2} (\psi +\eta)_{z\bar{z}}.
 \end{equation}
Moreover, since from \eqref{recgh}, \eqref{bjh2} we know that
$4|q|^2/\landa =\ro^2 |G_z|^2$, we infer from \eqref{bjh11} that
 \begin{equation}\nonumber
\psi=\psi +\eta +\frac{2}{\ro^2 |G_z|^2} (\psi +\eta)_{z\bar{z}}.
 \end{equation}
This lets us conclude by means of \eqref{bjh9}, \eqref{bjh105} that the
Bryant surface we started with is recovered in a neighbourhood of the
curve as
$\psi=F\Omega F^*:W\subseteq D \flecha \H^3$, where $F:W\subseteq
D\flecha \SL$ is the meromorphic curve in \eqref{bjh88}, and
$\Omega:W\subseteq D\flecha \He$ is the map (possibly with entries of infinite value at
some points)
 \begin{equation}\label{bjh13}
\Omega= \left(\def\arraystretch{2}\begin{array}{cc} \ro
+\displaystyle\frac{2\ro_{z\bar{z}}}{\ro^2 |G_z|^2} & \displaystyle
\frac{2\ro_z}{\ro^2 G_z}
\\ \displaystyle
\frac{2\ro_{\bar{z}}}{\ro^2 \overline{G_z}} & \displaystyle\frac{2}{\ro}
\end{array}\right).
 \end{equation}
Observe that ${\rm det} (\Omega)=1$, as $\ro$ satisfies the modified
Liouville equation.

To sum up, we have proved that the Bryant surface we started with is
completely determined in a neighbourhood of
$\beta$ by the Björling data $\beta,V$ and, furthermore, can be expressed in terms of
them. The analyticity of Bryant surfaces provides then uniqueness.

To proof existence we begin with Björling data $\beta,V$ and define
$\nu=\beta+V$, with values in the positive light cone $\N^3$. Again, we may assume
that $\nu_0(s)+\nu_3(s)\neq 0$ for all $s\in I$. Let $G:I\flecha \C$ be
given by
$$G(s)=\frac{\nu_1(s)-i\nu_2(s)}{\nu_0(s)+\nu_3(s)}, $$ take a meromorphic extension
$ G(z)$ of $G(s)$ and consider the only solution
$\ro:D\subseteq\C\flecha [0,+\8)$ to the Cauchy problem
\begin{equation}\nonumber\left\{\def\arraystretch{1.2} \begin{array}{lll}\displaystyle  4\left(\log
\ro\right)_{z\bar{z}}&=&-\ro^2 |G_z|^2,
\\ \ro(s,0)&=& v(s),\\
\ro_z(s,0)&=& w(s),
\end{array}\right.
\end{equation}
where $v(s),w(s)$ are given by \eqref{bjh3} and \eqref{bjh8},
respectively. Here $D\subseteq\C$ is an open subset containing
$I$, where $G(z)$ is also defined. Again, this Cauchy problem may be solved
via Corollary \ref{sollimo}, and its solution $\ro$ is recovered in terms
of $\beta,V$.

Let now $\cU\subseteq D$ be the open set
$\cU=D\setminus P$, with $$P=\{ \text{zeroes and poles of } G_z\} \cup
\{ \text{zeroes of } \ro \}.$$ Then we may define
$\psi=F\Omega F^*:\cU\flecha \H^3$, where $F,\Omega $ are given by
\eqref{bjh88} and \eqref{bjh13}, respectively. With this, define
$N:\cU\flecha \N^3$ as
 \begin{equation}\label{bjrec1}
N= F\left(\def\arraystretch{1.1}\begin{array}{cc} \ro & 0 \\ 0 & 0\\
\end{array}\right) F^*.
 \end{equation}
It comes plain that $N$ is conformal, and
$[N]:\cU\flecha \C_{\8}$ is meromorphic. In addition $\esiz N,\psi\esde=-1$,
and $\esiz N_z,\psi\esde=0$. Finally, since $[N]$ is meromorphic, we
conclude that $\psi$ must be a Bryant surface at its regular points, with
$N=\psi+\eta$.

We only have left to check that $\psi$ solves indeed the Cauchy problem
for the initial data
$\beta,V$. In this way we will also ensure that $\psi$ is regular in a neighbourhood
of $\beta$. First of all, from \eqref{bjh3} and \eqref{bjrec1} it follows
directly that
$N(s,0)=\beta(s)+V(s)$. So, we just need to verify
$\psi(s,0)=\beta(s)$.

The first step for this is to note that for arbitrary Björling data
$\beta,V$ the identity
 \begin{equation}\nonumber
(\nu_0+\nu_3)^2 G'\left\{\beta_1 +i\beta_2
 -\bar{G}(\beta_0+\beta_3)\right\} = \nu_0' +\nu_3' +i\left\{(V\ex
 \nu')_0 +(V\ex \nu')_3\right\},
 \end{equation}
holds. From this expression, as a consequence of \eqref{bjh8} we get
 \begin{equation}\label{bjrec2}
 2\ro_z (s,0)=(\nu_0+\nu_3)^2 G' \left\{\beta_1 +i\beta_2
 -\bar{G}(\beta_0+\beta_3)\right\}.
 \end{equation}
Let us write the matrix $\Omega$ in \eqref{bjh13} as
$$\Omega = \left(\begin{array}{cc} \omega_1 & \omega_2 \\
 \overline{\omega_2} & \omega_3\\
\end{array}\right).$$ Then \eqref{bjh3}, \eqref{bjrec2} lead to
 \begin{equation}\label{bjrec3}
\omega_2(s,0)= \frac{2\ro_z}{\ro^2 G'}(s,0) = \beta_1 +i\beta_2
-(\beta_0+\beta_3) \bar{G}.
 \end{equation}
Taking modulus we obtain
 \begin{equation*}
 |\omega_2|^2 (s,0) = \left|\beta_1 +i\beta_2
 -\bar{G}(\beta_0+\beta_3)\right|^2,
 \end{equation*}
and since ${\rm det}(\Omega)=1$, i.e. $2\omega_1 = \ro (1+|\omega_2|^2)$,
we have that
 \begin{equation*}
 2 \omega_1(s,0)=(\nu_0 +\nu_3) \left\{ 1+\left| \beta_1 +i\beta_2
 -\bar{G}(\beta_0+\beta_3)\right|^2\right\}.
 \end{equation*}
But additionally it can be checked that
 \begin{equation}\nonumber
 2(\beta_0+\beta_3) =(\nu_0 +\nu_3) \left\{ 1+\left| \beta_1 +i\beta_2
 -\bar{G}(\beta_0+\beta_3)\right|^2\right\},
 \end{equation}
and this tells that $ \omega_1(s,0)=\beta_0(s)+\beta_3(s)$. Hence, from
\eqref{bjrec3},
$$\omega_2(s,0) + \omega_1(s,0)\overline{G(s)}=\beta_1(s)+i\beta_2(s).$$ Thus, we obtain
that $\psi(s,0)=\beta(s)$. This ends up the proof.
\end{proof}

\begin{remark}
As Corollary \ref{sollimo} indicates, the solution of the Cauchy problem
\eqref{caulb} relies on solving the Cauchy problem for Liouville equation
\eqref{caulio} with $c=1$ and the initial data \eqref{insolli}, being
$f=G_z$.

A direct computation taking into account that $\esiz
(\psi+\eta)_z,(\psi+\eta)_{\bar{z}}\esde = \phi /2$ shows that such
initial data $a(s),b(s)$ are expressed in terms of the Björling data as
 \begin{equation}\label{datbjli}
 \def\arraystretch{1.2} \begin{array}{lll}
a&=& \esiz \nu',\nu' \esde,\\ b &=& \esiz \nu'',\nu'\esde + i\esiz
\nu',\nu'\esde \left( \displaystyle\frac{\big(V\wedge\nu'\big)_0 +
\big(V\wedge\nu'\big)_3}{\nu_0+\nu_3} + {\rm Im}
\left(\displaystyle\frac{G''}{G'}\right)\right).
\end{array}
 \end{equation}
Therefore, the Cauchy problem for Liouville equation that we have to solve
in order to obtain the only Bryant surface with Björling data
$\beta,V$ is \eqref{caulio} for $c=1$ with the initial conditions \eqref{datbjli}.
\end{remark}

It is possible to simplify the description of the only solution to the
Cauchy problem for Bryant surfaces with initial data $\beta,V$, by means
of two fundamental equations of the theory.

The first one is due to Umehara and Yamada \cite{UY1}. Let
$\psi:\Sigma\flecha \H^3$ be a Bryant surface with hyperbolic Gauss map $G$ and
Hopf differential $Q$, and denote by $ds^2,K$ its metric and curvature,
respectively. Then we may define on the universal cover
$\widetilde{\Sigma}$ of $\Sigma$ the \emph{secondary Gauss map} $g$
as the developing map of the curvature one pseudo-metric $-K ds^2$. In
other words, $g:\widetilde{\Sigma} \flecha \C\cup \{\8\} $ is defined by
$$-K ds^2 = \frac{4|dg|^2}{(1+|g|^2)^2}.$$ Then, the Umehara-Yamada
differential equation indicates that $G,Q,g$ are related on
$\widetilde{\Sigma}$ by
 \begin{equation}\label{uya}
\cS(g)-\cS(G) = -2Q.
 \end{equation}
Here, $\cS(g)=\{g,z\} dz^2$ and $\cS(G)=\{G,z\} dz^2$, where $z$ is an
arbitrary global complex parameter on $\widetilde{\Sigma}$.

The other fundamental equation we shall use is Small's formula \cite{Sma},
in the form exposed in \cite{GMM2}. This is nothing but an elaborated
version of the Bryant representation in \cite{Bry}. If
$G,g$ denote the hyperbolic and secondary Gauss maps of the Bryant surface
$\psi:\Sigma\flecha \H^3$, then $\psi =F
F^*:\widetilde{\Sigma}\flecha \H^3$, where
$F:\widetilde{\Sigma}\flecha \SL$ is the holomorphic curve
 \begin{equation}\label{small}
F= \left(\begin{array}{cc} dC/ dg & C - g\,dC/dg \\ dD/ dg & D - g \, dD/
dg\\
\end{array}\right), \hspace{.6cm} C=i \sqrt{dg / dG}, \hspace{0.2cm} D=i G\,\sqrt{dg / dG}.
 \end{equation}
Thus, $\psi$ is recovered explicitly in terms of $G,g$. With all of this,
we have as a consequence of Theorem \ref{solbjbr}:
 \begin{cor}\label{bjfaci}
Let $\beta,V$ be Björling data, define $\nu=\beta+V$ and
$G:I\flecha \C\cup \{\8\} $ by \eqref{genues}. Let $g:I\flecha \C\cup
\{\8\} $ be an arbitrary solution of $$\{g,s\} =\{G,s\} +\esiz
\nu'(s),\beta'(s)-i V(s)\wedge \beta'(s)\esde.$$ The only solution to the
Cauchy problem for Bryant surfaces with initial data $\beta,V$ is
constructed as $\psi=F F^*:\cU\subseteq\C\flecha \H^3$, where
$F:\cU\subseteq\C\flecha \SL$ is the holomorphic curve in \eqref{small}, being
$G(z),g(z)$ meromorphic extensions of $G,g$.
 \end{cor}

\begin{remark}
Corollary \ref{bjfaci} can also be proved directly by repeating some of
the computations in the proof of Theorem \ref{solbjbr}, but without the
need to consider Liouville equation. Nevertheless, there are at least
three basic reasons which justify the approach given in Theorem
\ref{solbjbr}. One, that the interaction between the two different
approaches provides a geometric resolution of the Cauchy problem for
Liouville equation with $c=1$ (see Section 5), which is more explicit than
the analytic one given in Section 2. Two, that the approach in Theorem
\ref{solbjbr} yields, among other results, an explicit construction of
Bryant surfaces out from one of its planar geodesics (see next Section).
And three, that the technique used in Theorem \ref{solbjbr} is very
flexible, and can be applied to solve the Cauchy problem for other classes
of surfaces associated to Liouville equation. This is the case, for
instance, of surfaces in $\H^3$ whose mean and Gaussian curvature $H,K$
verify a linear relation of the type $-2a H +b(K-1)=0$, with $|a+b|=1$,
$a,b$ constants (see \cite{GMM2}).
\end{remark}

The solution of the Cauchy problem for Bryant surfaces specified in
Theorem \ref{solbjbr} (and in Corollary \ref{bjfaci}) can be seen as a
complex representation for this type of surfaces, in which the input is a
pair of Björling data. It is interesting to remark that this
representation can be reformulated in a global fashion, in terms of
meromorphic data on a Riemann surface.

Specifically, let $\Sigma$ be a simply connected Riemann surface and
$\Gamma\subset \Sigma$ a regular analytic curve. Then we can define
\emph{Björling data} $\beta:\Gamma\flecha \H^3$ and $V:\Gamma\flecha
\S_1^3$ along $\Gamma$ such that $\esiz \beta,V\esde = \esiz
d\beta,V\esde=0$. From these data, we may consider $\nu=\beta
+V:\Gamma\flecha \N^3$, as well as $G:\Gamma\flecha \C\cup \{\8\} $ given
by $$G=\frac{\nu_1 -i \nu_2}{\nu_0 +\nu_3}.$$ We can also introduce a
complex $2$-form $Q$ along $\Gamma$ as $$Q=-\frac{1}{2} \esiz d \nu,
d\beta - i \beta \times \nu \times d\beta \esde.$$ Now, denote
$\cS(G)=\{G,s\} ds^2$, which is a complex $2$-form along $\Gamma$, and
consider a map
$g:\Gamma\flecha \C\cup \{\8\} $ such that $\cS(g)-\cS(G)=-2Q$ along
$\Gamma$. Suppose that:
 \begin{enumerate}
 \item
$g,G$ have meromorphic extensions to $\Sigma$, and $Q$ extends to $\Sigma$
as a global holomorphic $2$-form.
 \item
The poles of $g$ of order $k$ agree with the zeroes of $Q/dg$ of order
$2k$.
 \end{enumerate}
Then the map $\psi =F F^*:\Sigma\flecha \H^3$, where $F:\Sigma\flecha \SL$
is given by the Small-type formula \eqref{small}, is a Bryant surface with
$\psi (\Gamma)=\beta$ and $\eta(\Gamma)=V$. Here $\eta:\Sigma\flecha
\S_1^3$ is the unit normal to $\psi$ in $\H^3$. Of course, the converse
trivially holds.

The metric of the surface in this global conformal representation, as well
as the metric in the representation given by Theorem \ref{solbjbr}, has a
quite complicated expression in terms of the data $\beta,V$. Indeed, in
both cases this metric is given by $$ds^2 = \left(1+|g|^2\right)^2 \left|
Q/dg \right|^2,$$ and to find $g$ one needs to solve the Umehara-Yamada's
differential equation \eqref{uya}.

Nevertheless, it is known \cite{Yu} that the metric $ds^2$ of the Bryant
surface is complete (resp. non-degenerate) if and only if the \emph{dual
metric} $${ds^2}^{\sharp} = \left(1+|G|^2\right)^2 \left| Q/ dG\right|^2$$
is complete (resp. non-degenerate). Of course, this dual metric is much
simpler to handle through the Björling data $\beta,V$, because we have
explicit expressions which recover $G,Q$ in terms of
$\beta,V$ \eqref{himabj}, \eqref{bjh6}.

\section{Applications}

In this Section we view the solution to the Cauchy problem in Theorem
\ref{solbjbr} as a conformal representation for Bryant surfaces, and
establish several consequences regarding their geometry.

It is well known that if a Bryant surface meets a hyperbolic plane in
$\H^3$ orthogonally,  then it is symmetric with respect to that plane. First of
all, and as an immediate consequence of Theorem \ref{solbjbr}, we have the
following generalization of the above symmetry principle.
 \begin{teo}[Generalized symmetry principle]\label{simbry}
Any symmetry in the initial data of the Cauchy problem for Bryant surfaces
generates a global symmetry of the resulting Bryant surface.
 \end{teo}
\begin{proof}
Let $\Phi$ be a symmetry of the Björling data $\beta,V$, i.e. $\Phi$ is a
positive rigid motion of $\H^3$ and there is an analytic diffeomorphism
$\Psi:I\flecha I$ such that $\beta \circ\Psi = \Phi \circ \beta$ and
$V\circ \Psi = d\Phi \circ V$. If $\tilde{\Psi}$ is a holomorphic
extension of $\Psi$, and
$\psi$ is the only solution to the Cauchy problem for Bryant surfaces and the initial data $\beta,V$,
then the maps $\psi \circ \tilde{\Psi}$ and $\Phi \circ \psi$ are Bryant
surfaces with the same initial data. Thus $$\psi \circ \tilde{\Psi} = \Phi
\circ \psi,$$ i.e. $\Phi$ is a global symmetry of the Bryant surface
$\psi$.
\end{proof}

If in this symmetry principle we make the choices $\Phi={\rm Id}$ and
$\Psi(s)=s+T$ for some $T>0$, we obtain a general resolution of the period problem
for Bryant surfaces with the topology of a cylinder.
 \begin{cor}\label{perbry}
Let $\beta(s),V(s)$ be $T$-periodic Björling data. The Bryant surface that
solves the Cauchy problem for these initial data via Theorem \ref{solbjbr}
has the topology of a cylinder near $\beta$, and its fundamental group is
generated precisely by $\beta$.

Conversely, any Bryant surface with the topology of a cylinder is
recovered in this way.
 \end{cor}

This Corollary may be seen as a period problem-free holomorphic
representation for Bryant cylinders. As solving the period problem is the
fundamental and hardest step in classifying families of Bryant surfaces
with non-trivial topology, the above representation may be useful to
provide classification results in the case in which the surfaces have the
topology of a cylinder.

In this direction, we provide next a general description free of the
period problem for complete Bryant surfaces with finite dual total
curvature, genus zero and two ends.

\begin{cor}
Let $\beta(s),V(s)$ be $2\pi$-periodic Björling data satisfying
\begin{enumerate}
\item[$a)$]
$G(s)$ given by \eqref{genues} is a quotient of trigonometric polynomials.
 \item[$b)$]
$q(s)$ as in \eqref{cuese} as well as $h(s)=q(s)/G'(s)$ are trigonometric
polynomials.
 \item[$c)$]
The zeroes of the entire extension $h(z)$ of $h(s)$ are of order $2k$,
agreeing with the poles of order $k$ of $G(z)$, and at each end of the
Riemann surface $\C /2\pi \Z$ the (well-defined) extensions of
$h(z)$ and $h(z)G(z)^2 $ do not both vanish simultaneously.
\end{enumerate}
Then the Bryant surface that solves the Cauchy problem for the initial
data $\beta,V$ is complete, and has finite dual total curvature, genus
zero and two ends.

Conversely, any complete Bryant surface of finite dual total curvature
with genus zero and two ends is recovered in this way.
\end{cor}
\begin{proof}
Let $\beta,V$ be Björling data in the above conditions. Then, the tensor
$$ (1+|G(z)|^2)^2 |h(z)|^2 |dz|^2$$ is a regular Riemannian metric which
is well defined on the Riemann surface $\C /2\pi \Z$. Now, our analysis on
Section 4 assures that this metric is the dual metric of the solution
$\psi$ to the Cauchy problem for Bryant surfaces with initial data $\beta,V$.
In addition, since by Corollary \ref{perbry} this surface is homeomorphic
to a cylinder, we get that $\psi$ is parametrized as a conformal map from
$\C /2\pi \Z$ into $\H^3$, and by $c)$ is a complete regular Bryant surface. Now,
condition a) shows that the hyperbolic Gauss map $G(z)$ of $\psi$ extends
meromorphically to the compactification of $\C /2\pi \Z$. Thus the image
of $G:\C /2\pi \Z\flecha \C\cup \{\8\} \equiv \S_{\8}^2$ has finite area
(counted with multiplicities), i.e. $\psi$ has finite dual total
curvature.

Conversely, given a complete Bryant surface $\psi$ of finite dual total
curvature with genus zero and two ends, the Riemann surface in which it is
defined is a sphere with two points removed. Thus, we can assume that it
is parametrized in $\C /2\pi \Z $. Let us denote $z=s+it$, and consider
$\beta(s)=\psi (s+it_0)$ and $V(s)=\eta(s +it_0)$, being $t_0\in \R$
arbitrary and $\eta:\C /2\pi \Z \flecha \S_1^3$ the unit normal to
$\psi$. Then $\beta,V$ are $2 \pi$-periodic Björling data, and $\psi$ is
the solution to the Cauchy problem for Bryant surfaces and these initial
data. Finally, by the regularity, the completeness and the finite total
dual curvature condition of $\psi$, the data $\beta,V$ must satisfy $a)$,
$b)$ and $c)$.
\end{proof}
For the case of Bryant surfaces of finite total curvature, we get:
 \begin{cor}\label{ctfdu}
Let $\beta(s),V(s)$ be $2\pi$-periodic Björling data satisfying
 \begin{enumerate}
\item
$G(s)$ extends to a ($2\pi$-periodic) meromorphic function on $\C$.
\item
$q(s)$ defined in \eqref{cuese} is a trigonometric polynomial.
 \item
The metric $(1+|G(z)|^2)^2|q(z)/ G'(z)|^2 |dz|^2$, which is defined on
$\C / 2\pi \Z$ by its construction, is regular and complete.
 \item
If $\phi$ is the only solution to the Cauchy problem for Liouville
equation with initial data \eqref{datbjli} (which is globally defined on
$\C / 2\pi \Z$), there are positive numbers $\gamma_1,\gamma_2$ such that
the following two limits exist for all $s$:
$${\rm lim}_{t\flecha \8} \frac{\phi (s+it)}{\exp(2\gamma_1 t)},
\hspace{1cm} {\rm lim}_{t\flecha -\8} \frac{\phi (s+it)}{\exp(-2\gamma_2
t)} $$
 \end{enumerate}
Then the Bryant surface that solves the Cauchy problem for the initial
data $\beta,V$ is complete, and has finite total curvature, genus zero and
two ends.

Conversely, any complete Bryant surface of finite total curvature with
genus zero and two ends is recovered in this way.
 \end{cor}
\begin{proof}
The first three conditions together with the periodicity of the Björling
data indicate that the resulting Bryant surface is regular, complete, has
the topology of a cylinder and is parametrized on $\C /2\pi \Z$. The
fourth condition implies that $$\iint_{(0,2\pi )\times \R} \phi (s,t) ds
dt < +\8 .$$ That is, the surface has finite total curvature.

The converse follows the steps in the proof of Corollary \ref{ctfdu},
noting that a complete Bryant surface with finite total curvature and the
topology of a cylinder is conformally parametrized on a twice punctured
sphere. We remark that the fourth condition must hold on any complete
Bryant cylinder with finite total curvature $\psi:\C /2\pi \Z\flecha \H^3$
with respect to any Björling data of $\psi$ of the type $\beta(s)=\psi(s+i
t_0)$, $V(s)=\eta (s+it_0)$, since in that case the pseudo-metric
$\phi$ has conical singularities at the ends in $\C /2\pi \Z$ (see
\cite{Bry}).
\end{proof}

Another application of the generalized symmetry principle concerns singly
periodic Bryant surfaces. A surface in $\H^3$ is \emph{singly periodic}
provided it is invariant under the action of a cyclic group
$\cG$ of isometries of $\H^3$ that acts proper and discontinuously. In this case,
the surface may be regarded in the obvious way as immersed in the
hyperbolic $3$-manifold $\H^3 /\cG$.

For this, we start with Björling data $\beta,V$ for which there is an
isometry $\Phi$ of $\H^3$ such that: (a) $\Phi$ is a symmetry of
$\beta,V$, and (b) the action of the cyclic group $\cG$ generated by
$\Phi$ is proper and discontinuous.

Then it follows from the generalized symmetry principle that the solution
$\psi$ to the Cauchy problem for the initial data $\beta,V$ is singly
periodic, and can be regarded as an immersed surface in
$\H^3 /\cG$. Moreover, if $\pi:\H^3\flecha \H^3 / \cG$ is the canonical projection, then
$\pi \circ \psi$ has the topology of a cylinder in $\H^3 /\cG$, with
fundamental group generated by $\pi\circ \beta$.

Obviously, this process can be reversed: if $\psi:M^2\flecha \H^3$ is a
singly periodic Bryant surface, we may view it as
$\tilde{\psi}:M^2/\Lambda \flecha \H^3 /\cG$, where $\Lambda$ is a
cyclic group of isometries of $M^2$, induced by $\cG$ via $\psi$. If
$M^2 /\Lambda$ has the topology of a cylinder, there is some regular
curve $\beta$ on the surface in $\H^3$ such that $\pi\circ \beta$
generates the fundamental group of $M^2 /\Lambda$, and the above
construction recovers the immersion $\psi$.

The generalized symmetry principle together with Theorem \ref{solbjbr} can
be used to classify Bryant surfaces which are invariant under
$1$-parameter groups of rigid motions in $\H^3$, without the need to solve any
differential equation. To begin with, we consider Bryant surfaces which
are \emph{hyperbolic invariant}, i.e. they are invariant under the group
of hyperbolic translations along a geodesic in
$\H^3$. Any such group is, up to rigid motions, of the form $$A_{\alfa}=
\left(\def\arraystretch{1}
\begin{array}{cccc}
 \cosh \alfa & 0 & 0 & \sinh \alfa \\
 0 & 1 & 0 & 0 \\
 0 & 0 & 1 & 0 \\
\sinh \alfa & 0 & 0 & \cosh \alfa
\end{array}\right), \hspace{0.3cm} \alfa\in \R,$$ and its \emph{orbits} are
\emph{hyperbolic circles}, with the exceptional case of the geodesic that
defines the axis.

In addition, a direct computation shows that the unit normal in
$\H^3$ of a hyperbolic invariant surface along any of its orbits verifies
that its projection over the hyperbolic plane which is orthogonal to such
orbit is constant.

By means of this property and Theorem \ref{solbjbr}, the following example
classifies all hyperbolic invariant Bryant surfaces.

\begin{eje}\label{invhi}
\emph{Let $\beta(s)$ be a hyperbolic circle (or a geodesic) in
$\H^3$ which is an orbit of a hyperbolic invariant Bryant surface. Up to a rigid motion,
and due to the above property of hyperbolic invariant surfaces, the curve
$\beta(s)$ and the unit normal of the surface along
$\beta(s)$, denoted
$V(s)$, are} $$\beta(s)=\big(a \cosh s, b,0,a \sinh s\big),
\hspace{0.5cm} V(s)=(\landa \cosh s,c,d,\landa \sinh s),$$ \emph{where
$a,b,c,d,\landa$ verify
$a^2-b^2=1$, $a>0$, $-\landa^2 +c^2 +d^2 =1$, $a\landa = bc$.
Conversely, the generalized symmetry principle ensures that the solution
to a Cauchy problem of this type must always be a hyperbolic invariant
Bryant surface. Now, from Theorem \ref{solbjbr} we know that the
hyperbolic Gauss map and the Hopf differential of this Bryant surface are
} $$Q= -\frac{a}{2} (a+\landa) dz^2, \hspace{0.5cm} G(z)=k_1
e^{-z},\hspace{0.3cm} \text{\emph{with}} \hspace{0.3cm} k_1=\frac{b+c -i
d}{a+\landa}.$$ \emph{Once here, the Umehara-Yamada's relation \eqref{uya}
gives}$$g(z)=\exp \left(i \sqrt{2k_2}z\right) ,\hspace{0.5cm} Q=
-\frac{1}{2} \left(k_1 + \frac{1}{2}\right) dz^2,$$ \emph{with
$k_2= -1/2 + a(a+\landa)$ ($\neq 0$). Now we may recover in explicit
coordinates the immersion via Small's formula \eqref{small}. The dual
metric is}
$$ds^{2\sharp}= \left|\frac{2k_2+1}{4k_1}\right|^2 \left( 1 +  |k_1|^2
e^{-2s} \right)^2 e^{2s} |dz|^2,\hspace{0.5cm} z=s+it,$$ \emph{which is
clearly regular and complete ($k_2=-1/2$ gives the horosphere, as well as
$k_1=0$). Thus, we have obtained all hyperbolic invariant Bryant surfaces, which are
parametrized as maps from $\C$ into $\H^3$. All of them are regular,
complete and of infinite total curvature (except for horospheres).}
\end{eje}

 \begin{figure}[h]
  \begin{center}
    \begin{tabular}{cc}
\includegraphics[clip,width=7cm]{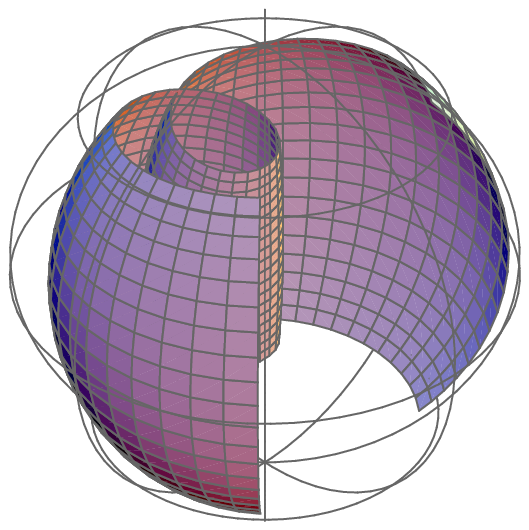} &
\includegraphics[clip,width=7cm]{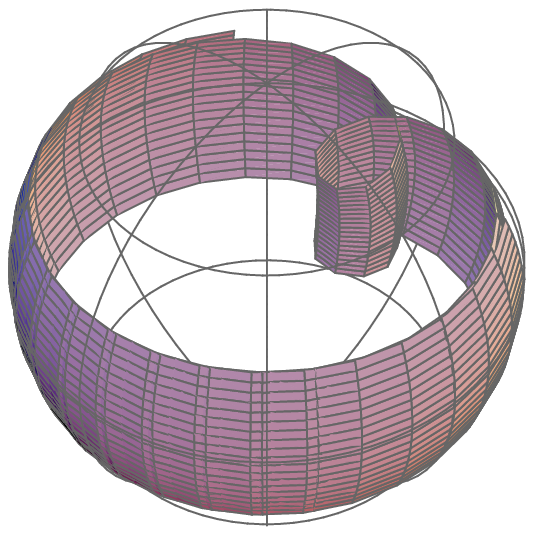}
\end{tabular}
 \end{center}
\caption{Hyperbolic invariant Bryant surfaces in the Poincaré model
containing and not containing the axis of the hyperbolic translation}
\end{figure}

 \begin{eje}[Catenoid cousins]\label{caco}
\emph{Catenoid cousins are the fundamental examples in the theory of
Bryant surfaces. Here, we will show that they admit a simple construction
in terms of Björling data. Indeed, next we find via Theorem \ref{solbjbr}
the Bryant surfaces that have a circle as a planar geodesic. We remark
that, depending on the radius of the circle, the resulting surfaces will
have distinct geometries. This contrasts with the case of minimal surfaces
in $\R^3$.}

\emph{We begin with a circle in $\H^3$,}
$$\beta(s)=\left( c, b \cos s,b \sin s, 0\right), \hspace{0.5cm} c=\sqrt{1+b^2}, \
b>0.$$\emph{If a Bryant surface contains
$\beta(s)$ as planar geodesic, its unit normal in $\H^3$ along
$\beta(s)$ is}
$$V(s)=\varepsilon \left(b, c \cos s, c  \sin s, 0\right), \hspace{0.5cm}
\varepsilon =\pm 1.$$ \emph{Hence, its hyperbolic Gauss map and its Hopf
differential are} $$G(z)= \varepsilon e^{-iz}, \hspace{0.5cm} Q=
-\frac{b}{2}\left( b +\varepsilon \sqrt{1+b^2}\right) dz^2.$$ \emph{From
here, and again by Umehara-Yamada's differential equation, we get }
$$g(z)= \exp \left( i\sqrt{2k} z\right), \hspace{0.5cm} Q = -\frac{1}{2}
\left(k-\frac{1}{2}\right) dz^2,$$ \emph{being $k= 1/2 + b\left( b
+\varepsilon \sqrt{1+b^2}\right)$. The surface is recovered in explicit
coordinates via \eqref{small}.}

\emph{From the generalized symmetry principle and Corollary \ref{perbry}
these are rotation Bryant surfaces with the topology of a cylinder,
parametrized as $\psi:\C /2\pi \Z\flecha \H^3$. The secondary Gauss map
$g$ is in general multivalued on $\C /2\pi \Z$, and the dual metric
$ds^{2\sharp}$ is the metric of a catenoid in $\R^3$. Thus, we have
obtained the \emph{catenoid cousins} \cite{Bry,UY1}. They are complete
Bryant surfaces, with dual total curvature $-4\pi$, finite total curvature
$-4\pi \sqrt{2k}$ (admitting all possible negative values), with genus
zero and two ends. Depending on the radius of the circle we started with,
the resulting catenoid cousin is embedded or not.}
 \end{eje}
 \begin{figure}[h]
  \begin{center}
    \begin{tabular}{cc}
\includegraphics[clip,width=7cm]{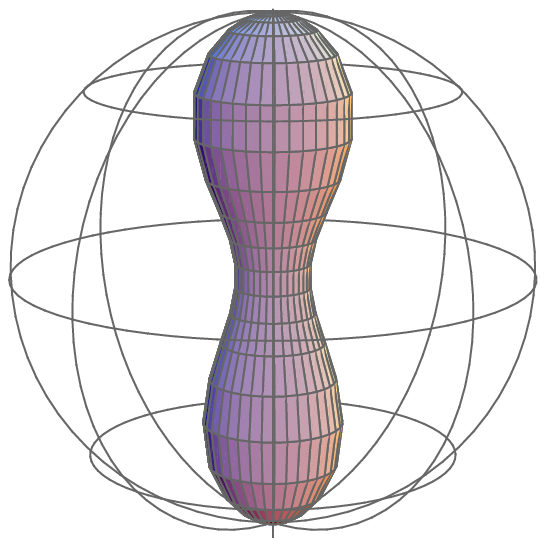} &
\includegraphics[clip,width=7cm]{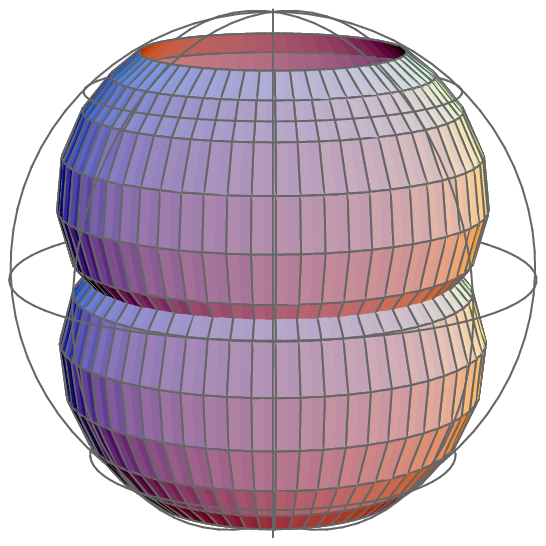}
\end{tabular}
 \end{center}
\caption{Embedded and non-embedded catenoid cousins in the Poincaré model}
\end{figure}

The techniques we are exploiting here can be use to provide a simple proof
of the fact that \emph{catenoid cousins are the only rotation Bryant
surfaces}. Nevertheless, we shall give instead the classification of
helicoidal Bryant surfaces, which of course will contain as a particular
case the above mentioned result.

A \emph{helicoidal motion} in $\H^3$ is a rotation composed with a
hyperbolic translation along the axis of the rotation. They form
continuous $1$-parametric groups sharing the axis. Specifically, any
\emph{helicoidal group} is given in adequate coordinates as
 \begin{equation}\nonumber
A_{\alfa}= \left(\def\arraystretch{1} \begin{array}{cccc}
 \cosh (\alfa s) & 0 & 0 & \sinh (\alfa s)\\
 0 & \cos s & -\sin s & 0 \\
 0 & \sin s & \cos s & 0 \\
\sinh (\alfa s)& 0 & 0 & \cosh (\alfa s)
\end{array}\right), \hspace{0.3cm} s\in \R
 \end{equation}
for some $\alfa\in \R$, called the \emph{angular pitch}. Of course,
helicoidal motions with $\alfa=0$ are precisely rotations. A \emph{helix}
in $\H^3$ is a non-geodesic orbit of a helicoidal motion. They are
geodesics of hyperbolic cylinders whose axis is that of the helicoidal
motion. A surface $S$ in $\H^3$ is \emph{helicoidal} if its image is
invariant under a helicoidal motion group. In that case $S$ is foliated by
helices, and meets at a constant angle the hyperbolic cylinders in which
those helices lie.

Putting together this fact and the solution to the Cauchy problem for
Bryant surfaces, we find the classification of helicoidal Bryant surfaces.

 \begin{teo}\label{helibry}
A Bryant surface is helicoidal with axis $\ell$ if and only if it meets a
hyperbolic cylinder $C$ of axis $\ell$ at a constant angle along a helix
of $C$.

Furthermore, this initial value problem in $\H^3$ is explicitly solved via
Corollary \ref{bjfaci}, and its solution has hyperbolic and secondary
Gauss maps and Hopf differential given by
 \begin{equation}\nonumber
G(z)= c_1 e^{-(\alfa + i)z}, \hspace{0.5cm} g(z)= e^{c_2 z},
\hspace{0.5cm} Q= c_3 dz^2
 \end{equation}
for constants $c_1,c_2,c_3\in \C\setminus\{0\}$. Particularly, any
helicoidal Bryant surface lies in the associate family of a catenoid
cousin.
 \end{teo}
\begin{proof}
The direct part is true for any helicoidal surface in
$\H^3$, while the converse follows from the generalized symmetry
principle.

Let us solve the above initial value problem in $\H^3$. For this, we start
with a helix, that we assume to be
 \begin{equation}\nonumber
 \beta(s)=\left( c \cosh(\alfa s), b\cos s, b\sin s, c\sinh (\alfa
 s)\right), \hspace{0.5cm} c^2-b^2=1, c>0.
 \end{equation}
This helix lies in the hyperbolic cylinder $C$ of equation
$-x_0^2+x_3^2 =-c^2$, and the unit normal of $C$ along $\beta(s)$ is
 \begin{equation}\nonumber
 \xi(s)= \left(b\cosh(\alfa s), c\cos s, c\sin s, b\sinh (\alfa
 s)\right).
 \end{equation}
From our hypothesis, the unit normal to the Bryant surface along
$\beta(s)$ is
 \begin{equation}\label{eluve}
V(s)= \cos\varphi \, \xi(s) + \frac{\sin \varphi}{|\beta'(s)|}
\beta(s)\times \beta'(s)\times \xi(s), \hspace{0.5cm} \varphi \in \R.
 \end{equation}
With this, we obtain directly from \eqref{himabj} that
$G(z)= c_1 e^{-(\alfa +i)z}$ for $$c_1= \frac{b +c\cos \varphi +i
\alfa \sin \varphi /\sqrt{c^2\alfa^2 +b^2}}{c +b\cos\varphi +b\sin \varphi
/\sqrt{c^2 \alfa^2 + b^2}} \in \C\setminus \{0\}.$$ The Hopf differential
is $Q= q(z) dz^2$ being $q(z)$ as in \eqref{jof}. We claim that $q(z)$ is
constant. Indeed, since
$\beta(s)$ is geodesic of $C$ with constant speed, $\beta''(s)$ lies in the plane spanned by
$\beta(s)$ and $\xi(s)$, and those span the normal plane
to $C$ in $\L^4$ along $\beta(s)$. In addition $\xi''(s)$ is also normal.
Hence, \eqref{eluve} yields
$${\rm det} (\beta,V,\beta', V')=\cos^2 \varphi \,{\rm det}
(\beta,\xi,\beta',\xi') +\frac{\sin^2\varphi}{|\beta'(s)|^2}\, {\rm det}
(\beta, \beta'\ex \xi, \beta', \beta'\wedge \xi'),$$ and we see that
$$\frac{d}{ds} \big( \esiz \beta'+V',\beta'\esde - i \,{\rm det}
(\beta,V,\beta',V')\big) =0.$$ Thus, $q(z)$ is constant and
$Q= c_3 dz^2$ for $c_3\in \C\setminus \{0\}$. Finally, from
Umehara-Yamada's relation \eqref{uya} we find that $g(z)=\exp (c_2 z) $
with $c_2\in \C\setminus\{0\}$. In particular we recover the immersion in
explicit coordinates by means of Small's formula \eqref{small}.

We have only left to check that all these surfaces are associated to some
catenoid cousin. To do so, we will show that any simply connected Bryant
surface such that $Q= a_1\, dw^2$ and
$\{g,w\} =a_2$ for some global conformal parameter $w$ and constants
$a_1,a_2\in \C\setminus \{0\}$ is associated to some catenoid cousin.

We already saw in Example \ref{caco} that any catenoid cousin possesses a
global conformal parameter $z$ such that its Hopf differential $Q_c$ and
its secondary Gauss map
$g_c$ satisfy
$$Q_c= -\frac{1}{2}\left( k-\frac{1}{2}\right) dz^2, \hspace{1cm}
\{g_c,z\} =k$$ being $k=1/2 +b(b +\varepsilon \sqrt{1+b^2})$,
$b>0,\varepsilon = \pm 1$. It is easy to see that $k>0$. In addition, when
$b\rightarrow \8$ if follows that $k\rightarrow 0$ if $\varepsilon=-1$,
and $k\rightarrow \8$ if $\varepsilon=1$. This guaranties that
$$\frac{q_c (z)}{\{g_c,z\}} = \frac{1}{2}\left( \frac{1}{2k}
-1\right)$$ assumes all positive values. Therefore we may select a
catenoid cousin with associated $k>0$ verifying
$$\frac{1}{4}\left( \frac{1}{2k}
-1\right)^2=\left|\frac{a_1}{a_2}\right|^2.$$ Consider besides the change
of parameter $w=\landa z$, with
$\landa^2 = k/a_2$. Then
$\{g_c,w\} =a_2$ and $$\frac{q_c (w)}{\{g_c,w\}} = \frac{1}{2}\left(
\frac{1}{2k} -1\right)=e^{i\theta} \frac{a_1}{a_2}$$ for some
$\theta\in\R$. Consequently, the Hopf differential of the catenoid
cousin is
$Q_c= e^{i\theta} a_1 \, dw^2$. Thus, the Bryant surface with
$Q= a_1\, dw^2$ and $\{g,w\}=a_2$ lies in the associate family of the catenoid
cousin that we have determined. This completes the proof.
\end{proof}

Catenoid cousins are not the only Bryant surfaces that can be explicitly
recovered by means of one of its planar geodesics. Indeed, as a
consequence of Corollary \ref{corgp} we obtain a general procedure to
construct in explicit coordinates Bryant surfaces containing a prescribed
curve as a planar geodesic.

\begin{teo}\label{geplabr}
Let $\beta(s)$ be a regular analytic curve lying in the hyperbolic plane
$\H^3\cap P$, with $P\equiv x_2=0$, let
${\bf e}=(0,0,1,0)$ and denote
 \begin{equation}\label{lanugp}
 \nu(s)=\beta(s) + \frac{\varepsilon}{|\beta'(s)|}\, \beta(s)\times
 \beta'(s)\times {\bf e}, \hspace{0.5cm} \varepsilon = \pm 1.
 \end{equation}
If $\beta$ is not a geodesic of $\H^3$, there are exactly two Bryant
surfaces that contain $\beta(s)$ as planar geodesic, and both of them are
constructed explicitly via Small's formula \eqref{small} for
 \begin{equation}\label{apsgs}
 G(z)= \frac{\nu_1(z)}{\nu_0(z)+\nu_3(z)}, \hspace{0.5cm} g(z)=\exp \left(
 i\int_{s_0}^z \esiz \nu'(\zeta),\nu'(\zeta)\esde^{1/2} d\zeta\right).
 \end{equation}
Here $\nu(z)$ is a holomorphic extension of $\nu(s)$, and $s_0\in I$ is
fixed and arbitrary.
\end{teo}
\begin{proof}
Given a Bryant surface that contains $\beta(s)$ as a planar geodesic, its
unit normal along $\beta(s)$, denoted $V(s)$, lies in
$\S_1^3\cap P$ and verifies $\esiz \beta,V\esde=\esiz
\beta',V\esde=0$. This indicates that
$\nu(s)=\beta(s) +V(s)$ is given by \eqref{lanugp}. Now, from Theorem
\ref{solbjbr} the hyperbolic Gauss map $G(z)$ is the one specified in
\eqref{apsgs}. Particularly, $G(s)\in \R$ for all $s$.

On the other hand, as $\beta(s),\nu(s)$ lie in $P$ we find that
$V\wedge \nu'=\beta\times\nu\times \nu'$ is collinear with $e$. Hence $\big(V
\wedge \nu'\big)_0 =\big(V \wedge \nu'\big)_3=0$, and the initial
conditions in \eqref{datbjli} satisfy
$2b(s)= a'(s)$. That is, we are left with the Cauchy problem
\begin{equation*}\left\{\def\arraystretch{1.2} \begin{array}{lll}\displaystyle  \lap \log
\phi &=&\displaystyle -2\phi, \\ \phi(s,0)&=& \esiz \nu'(s),\nu'(s)\esde,
\\ \phi_t(s,0)&=& 0 .
\end{array}\right.
\end{equation*}
Now, this Cauchy problem can be explicitly solved via Corollary
\ref{corgp}, and we obtain
 \begin{equation*}
\phi(s,t)=\frac{4|g_z|^2}{(1+|g|^2)^2} \hspace{0.5cm} \text{for}
\hspace{0.5cm} g(z)=\exp\left( i\int_{s_0}^z \esiz
\nu'(\zeta),\nu'(\zeta)\esde^{1/2} d\zeta\right).
 \end{equation*}
This ensures that the secondary Gauss map of the Bryant surface is exactly
the one given in \eqref{apsgs}. Finally, this surface is recovered via
\eqref{small}. The two choices of the sign
$\varepsilon=\pm 1$ determine the two unique Bryant surfaces containing
$\beta(s)$ as planar geodesic.
\end{proof}

Let us observe that if a Bryant surface meets a hyperbolic plane
orthogonally along a geodesic of $\H^3$, then its Björling data are, up to
a rigid motion, $$\beta(s)= (\cosh s, 0,0,\sinh s),\hspace{1cm}V(s)=(0,\pm
1,0,0).$$ That is, we obtain the only (up to congruences) hyperbolic
invariant Bryant surface in $\H^3$ that contains the axis of the
translation (see Example \ref{invhi}).

\begin{remark}
If $\beta$ is a planar geodesic of a Bryant surface, then this surface is
symmetric with respect to the hyperbolic plane in which
$\beta$ lies. This fact has been used in the theory to construct
examples, and shows the importance of planar geodesics on Bryant surfaces
\cite{Bob,Kar,RUY1}. For instance, it is pointed out in \cite{Bob} that
Bryant trinoids in $\H^3$ contain planar geodesics, and that the important
(and open) question of classifying embedded trinoids can be reduced to
determining when these planar geodesics are embedded.
\end{remark}

For the general case of Bryant surfaces containing a given curve as a
geodesic, we get the following result, whose proof follows from Theorem
\ref{solbjbr} after a direct computation.

 \begin{cor}
Let $\beta(s)$ be a regular analytic curve in $\H^3$ that is not a
geodesic. There exist exactly two Bryant surfaces that contain $\beta$ as
a pregeodesic. Both of them are constructed by solving the Cauchy problem
for Bryant surfaces with initial data $\beta$ and $$V=\pm \left( \beta'' -
\frac{\esiz \beta'',\beta'\esde}{\esiz \beta',\beta'\esde }\beta' - \esiz
\beta',\beta'\esde \beta \right) \left/ \left|\left| \beta'' - \frac{\esiz
\beta'',\beta'\esde}{\esiz \beta',\beta'\esde }\beta' - \esiz
\beta',\beta'\esde \beta \right|\right|\right. $$ by means of Theorem
\ref{solbjbr} or Corollary \ref{bjfaci}.
 \end{cor}

For the next application, we recall the \emph{Bryant representation}
\cite{Bry}, asserting that for any Bryant surface $\psi:\Sigma\flecha\H^3$
there exists a holomorphic null curve $F:\widetilde{\Sigma}\flecha \SL$ on
the universal cover $\widetilde{\Sigma}$ of $\Sigma$ such that the
identity $\psi=F F^*:\widetilde{\Sigma}\flecha \H^3$ holds. The matrix $F$
is in general multivalued on $\Sigma$, even though $F F^*$ is not. We say
that
$\psi$ \emph{lifts to a null curve} in $\SL$ if $F$ is well defined on
$\Sigma$. We remark that
$\psi$ lifts to a null curve in $\SL$ if and only if the secondary Gauss
map
$g:\widetilde{\Sigma}\flecha \C\cup\{\8\} $ is single valued on $\Sigma$
(see \cite{UYMath}). This property is important in connection with the
study of rigidity of Bryant surfaces. We say that a Bryant surface
$\psi:\Sigma\flecha \H^3$ is \emph{rigid} if any other mean curvature one immersion
from the Riemannian surface $\Sigma$ into $\H^3$ only differs from $\psi$
by a rigid motion. It was then proved in \cite{UYMath} that if $\psi$ does
not lift to a null curve in $\SL$, then it is rigid. In contrast, there
are rigid Bryant surfaces which lift to a null curve in $\SL$.

Let us also introduce the following definition:
 \begin{defi}
A pair of Björling data $\beta,V$ is \emph{admissible} provided the
function given in terms of $\beta,V$ by
 \begin{equation}\label{curad}
 \kappa (s)= \frac{-1}{\esiz \nu',\nu'\esde^{1/2}} \left(
\displaystyle\frac{\big(V\wedge\nu'\big)_0 +
\big(V\wedge\nu'\big)_3}{\nu_0+\nu_3} + {\rm Im}
\left(\displaystyle\frac{G''}{G'}\right)\right)
 \end{equation}
is well defined at all points, and non-constant.
 \end{defi}

Let $\psi:\Sigma\flecha \H^3$ be a Bryant surface with unit normal
$\eta$, and let $\gamma_1(s),\cdots,\gamma_n(s)$ be periodic curves in $\Sigma$
that generate the first homology group $H_1(\Sigma,\Z)$. We will assume
that the pair of Björling data
$\beta_j(s)=\psi(\gamma_j(s)),V_j(s)=\eta(\gamma_j(s))$ are admissible.
Let $\alfa_j(s):\R\flecha \S^2$ be the only curve in $\S^2$ with arclength
parameter $u_j(s)= \int^s \esiz \nu_j'(r),\nu_j'(r) \esde^{1/2} dr$ and
geodesic curvature $$\kappa_j(s)= \frac{-1}{\esiz
\nu_j',\nu_j'\esde^{1/2}} \left(
\displaystyle\frac{\big(V_j\wedge\nu_j'\big)_0 +
\big(V_j\wedge\nu_j'\big)_3}{(\nu_j)_0+(\nu_j)_3} + {\rm Im}
\left(\displaystyle\frac{G_j''}{G_j'}\right)\right)$$ Then we have

\begin{cor}\label{fole}
There is a finite-folded covering of $\psi$ that lifts to a null curve in
$\SL$ if and only if the curves
$\alfa_j(s)$ in $\S^2$ are all closed.
\end{cor}
\begin{proof}
First, observe that the derivative of the arclength parameter $u_j'(s)$ of
$\alfa_j(s)$, as well as $\kappa_j(s)$, are both periodic, with the period of
$\gamma_j(s)$. Let $T_j$ denote this period. Then, since $\kappa_j(s)$ is
not constant, it is easy to realize that $\alfa_j$ is closed if and only
if $\alfa_j(s)$ is periodic with period $l_j T_j$ for some $l_j\in \N$.

In addition, let $g:\widetilde{\Sigma}\flecha \C\cup \{\8\} $ be the
secondary Gauss map of $\psi$, and denote $\tilde{g}_j(s) =
g(\tilde{\gamma}_j(s))$. Here $\tilde{\gamma}_j(s):\R\flecha
\widetilde{\Sigma}$ is the lift to $\widetilde{\Sigma}$ of $\gamma_j(s)$.
Hence, from Theorem \ref{curvas} we know that $\tilde{g}_j(s)=\pi
(\alfa_j(s))$, $\pi$ denoting stereographic projection. Thus, $\alfa_j$ is
closed if and only if there is some $l_j\in \N$ such that
$\tilde{g}_j(s)$ is $l_j T_j$-periodic. This indicates that all curves
$\alfa_j$ are closed if and only if there exist $l_1,\dots,l_n \in \N$ so
that $g$ is single valued over the finite-folded covering of $\Sigma$ with
first homology group generated by $l_1\gamma_1, \dots, l_n \gamma_n $.
Thus, all $\alfa_j$ are closed if and only if a finite-folded covering of
$\psi$ constructed as above lifts to a null curve in $\SL$.
\end{proof}

The previous Corollary shows that the question of determining when does a
Bryant surface lift to a null curve in $\SL$ is strongly related to the
following problem posed by S.S. Chern (see \cite{GaMi1}): \emph{when is a
curve with periodic curvatures periodic?} In \cite{GaMi1} the authors used
a perturbation method to obtain several results regarding this problem for
the case of curves in $\S^2$ and, as a consequence, described the first
family of flat tori in $\R^4$ that comes out since the 19th century. Next,
we apply a modification of this method to Bryant surfaces.

Let $\cC$ denote the class of closed admissible Björling data $\beta,V$.
We shall use the parameter $s=u/m\ell$, where $u$ is the arclength
parameter of $\beta$, $\ell$ is the length of $\beta$, and $m$ is the
number of times that we have to trace $\beta$ so that $V$ also closes.
Then, we view these Björling data as a pair of maps
$(\beta(s),V(s)):[0,1]\flecha \H^3\times \S_1^3$. A topology on $\cC$ may
be defined through the norm
$$|| \cdot ||_{\cC} = ||\beta||_{\8} +||\beta'||_{\8}+||\beta''||_{\8}+||
V ||_{\8}+|| V'||_{\8}+||V''||_{\8}.$$ With this, we get

 \begin{cor}\label{ever}
There exists a continuous functional $\cA$ from $\cC$ into $\R$ such that
$\cA (\beta,V) \in \Q$ if and only if the Bryant cylinder that $\beta,V$ generate
has a finite-folded covering that lifts to a null curve in $\SL$.
 \end{cor}
\begin{proof}
Let $(\beta(s),V(s))\in \cC$. Then we can define in terms of them the maps
$\nu (s)=\beta(s)+V(s)$ and $G(s)$ as in \eqref{genues}. Once here, we may
assign to $\beta,V$ the pair of functions defined on $[0,1]$ given by
$$u(s)=\int^s \esiz \nu'(r),\nu'(r)\esde dr, \hspace{0.5cm} \kappa (s)=
\frac{-1}{\esiz \nu',\nu'\esde^{1/2}} \left(
\displaystyle\frac{\big(V\wedge\nu'\big)_0 +
\big(V\wedge\nu'\big)_3}{\nu_0+\nu_3} + {\rm Im}
\left(\displaystyle\frac{G''}{G'}\right)\right) ,$$ and then assign to
this new pair the only (up to congruences) curve $\alfa(s)$ in $\S^2$ with
arclength parameter $u(s)$ and geodesic curvature $\kappa(s)$. It is not
difficult to adapt the arguments in the proof of Theorem $18$ of
\cite{GaMi1} to ensure that these mappings are continuous.

Finally, we extend $\alfa(s):[0,1]\flecha \S^2$ to a curve defined in the
whole $\R$, and consider the only rigid motion $A_{\alfa}$ in $\S^2$ such
that $A_{\alfa} (\alfa(0))=\alfa(1)$, $A_{\alfa} (\alfa'(0))=\alfa'(1)$
and $A_{\alfa} (\alfa(0)\times \alfa'(0))=\alfa(1)\times \alfa'(1)$. Then
$\alfa(s+1)=A_{\alfa} (\alfa(s))$ holds for all $s\in \R$, since both curves
have the same geodesic curvature, and the same initial conditions. Observe
that $A_{\alfa}\in {\rm SO}(3)$. If $\theta_{\alfa}$ denotes the angle of
the rotation $A_{\alfa}$, it then comes clear that $\theta_{\alfa}/\pi \in
\Q $ if and only if $A_{\alfa}^q ={\rm Id} $ for some $q\in \N$, if and
only if $\alfa(s+ q)=\alfa(s)$ for some $q\in \N$ and for all $s\in \R$.
In addition, the mapping
$$\alfa(s)\mapsto \theta_{\alfa} \in \R$$ is continuous (see
\cite{GaMi1}). By putting together all of this, we get a continuous map
$\cA:\cC\flecha \R$ with the property that
$\cA(\beta,V)\in \Q$ if and only if the curve $\alfa$ constructed by means of $\beta,V$
is closed (recall that $\beta,V$ are admissible). The proof concludes then
by Corollary \ref{fole}

\end{proof}

This Corollary provides a perturbation procedure in the class of Bryant
cylinders that lift to a null curve in $\SL$, so that a prescribed curve
on the surface remains invariant in the process. Specifically, let $S$ be
a Bryant cylinder that lifts to a null curve in $\SL$, and let $\beta(s)$
be a regular analytic curve in $\H^3$ that generates the fundamental group
of $S$, and such that the Björling data $\beta,V$ of $S$ are admissible.
It follows then by Corollary \ref{ever} that $\cA(\beta,V)\in \Q$. It is
not difficult to convince oneself that the pair $\beta,V$ is not isolated
in $\cC$, and that one can perturb this pair within its connected
component in $\cC$ so that $\beta$ remains fixed, and $V(s)$ changes. By
continuity of the functional $\cA$ we obtain the existence of new Björling
pairs, which generate Bryant cylinders distinct from $S$ that pass through
$\beta$, and whose value under $\cA$ is rational. That is, we obtain a
process to deform $S$ that has the required properties.

Up to now we have avoided the situation in which the function
\eqref{curad} given in terms of $\beta,V$ is constant, since this case
behaves in a different way. However, as we saw in Theorem \ref{geplabr}
the function \eqref{curad} vanishes when we deal with planar geodesics of
Bryant surfaces. As this situation is of a special interest, next we
investigate the topic in Corollary \ref{ever} for the case of planar
geodesics.

\begin{cor}
Let $\psi:\Sigma\flecha \H^3$ be a Bryant surface, and suppose there exist
periodic curves $\gamma_1(s),\cdots ,\gamma_n(s)$ in $\Sigma$ with periods
$T_1,\cdots ,T_n$ that generate the first homology group $H_1(\Sigma,\Z)$, and such
that $\beta_j(s)=\psi (\gamma_j(s))$ are planar geodesics. The following
two conditions are equivalent:
 \begin{enumerate}
\item[i)]
$\psi$ lifts to a null curve in $\SL$.
\item[ii)]
for every $j=1,\cdots, n$ it holds $$\int_0^{T_j} \esiz \nu_j'(s),
\nu_j'(s)\esde^{1/2} ds \in 2\pi \Z,$$ where $\nu_j(s)$ is defined in
terms of $\beta_j(s)$ as \eqref{lanugp}.
 \end{enumerate}
\end{cor}
\begin{proof}
Let $g$ be the secondary Gauss map of $\psi$, and denote $\tilde{g}_j(s)=
g(\widetilde{\gamma}_j(s)):\R\flecha \C\cup \{\8\} $, where
$\widetilde{\gamma}_j(s)$ is the lift to $\widetilde{\Sigma}$ of $\gamma_j$. Then $g$ is
single valued on $\Sigma$ if and only if $\tilde{g}_j(s)$ is
$T_j$-periodic for all $j$. But now we know from Theorem \ref{geplabr}
that $$ \tilde{g}_j (s)=\exp \left(
 i\int_{s_0}^s \esiz \nu_j'(r),\nu_j'(r)\esde^{1/2} dr \right).$$
Therefore, $\tilde{g}_j(s)$ is $T_j$-periodic if and only if the condition
ii) holds. This ends up the proof.
\end{proof}

\section{Liouville equation revisited}

In this Section we provide a geometric resolution of the Cauchy problem
 \begin{equation}\label{caulior1}
\left\{\def\arraystretch{1.2} \begin{array}{lll}\displaystyle \lap (\log
\phi)&=&\displaystyle -2 \phi, \\ \phi(s,0)&=& a(s),\\ \phi_t(s,0)&=&
d(s),
\end{array}\right.
 \end{equation}
by means of the two alternative resolutions of the Cauchy problem for
Bryant surfaces exposed in Section 3. In this resolution, and opposite to
the analytic approach developed in Section 2, we do not need to complexify
Liouville equation. In addition, we recover the \emph{developing map}
$g$ that describes the solution $\phi$ via formula \eqref{sollio} for $c=1$.

To begin with, let $a(s),d(s):I\flecha \R$ be the (analytic) initial data
of the Cauchy problem \eqref{caulior1}. As usual, we assume that
$a(s)$ is positive. Let $\nu(s):I\flecha \N^3$ be a regular analytic
curve in the positive light cone, parametrized so that $\esiz
\nu'(s),\nu'(s)\esde = a(s)$, and which verifies $\nu_0(s) +\nu_3(s) \neq
0$ for all $s\in I$. Next, we define $\cL (s):I\flecha \R$ as
 \begin{equation}\nonumber
 \cL(s)= -\left(\nu_0(s) + \nu_3(s)\right) \left(
 \displaystyle\frac{d(s)}{2 a(s)} + {\rm Im} \left(
 \displaystyle\frac{G''(s)}{G'(s)}\right) \right),
 \end{equation}
and let $\cF(s):I\flecha \L^4$ be an analytic curve defined by equations
 \begin{equation}\label{condefe}
\esiz \cF(s),\nu(s)\esde = \esiz \cF(s),\nu'(s)\esde =0, \hspace{0.5cm}
\esiz \cF(s),\cF(s)\esde =a(s), \hspace{0.5cm} \cF_0(s) +\cF_3(s) =
\cL(s).
 \end{equation}
The first three conditions ensure that, for each particular $s$, $\cF(s)$
lies in the intersection of a $2$-dimensional de Sitter space centred at
the origin of $\L^4$, with a degenerate vector hyperplane of $\L^4$. This
intersection is made up by a pair of straight lines which lie in the
degenerate vector hyperplane of equation $\esiz \nu(s),W\esde =0$. But
now, since $\nu_0(s) +\nu_3(s) \neq 0$, this hyperplane is not parallel to
the degenerate affine hyperplane $\cF_0 (s) +\cF_3(s) = \cL (s)$.
Consequently, for each $s$ we obtain two possible values of $\cF(s)$
satisfying \eqref{condefe}, one for each straight line. This provides the
existence of such $\cF(s)$, and shows how to construct it by purely
algebraic manipulations.

Let us finally define $\alfa(s)$ as the only curve $\alfa(s):I\flecha
\N^3$ with $\alfa(s)\times \nu(s)\times \nu'(s) =\cF(s)$ and $\esiz
\alfa(s),\nu'(s)\esde =0$. Then we have
\begin{teo}
Let $a(s),d(s):I\flecha \R$ be real analytic functions, with $a(s)$
positive. The only solution to the Cauchy problem \eqref{caulior1} is
 \begin{equation}\label{resceuno}
 \phi(s,t)= \frac{4 |g'(z)|^2}{(1+|g(z)|^2)^2}, \hspace{1cm} z=s+it,
 \end{equation}
where $g(z)$ is a meromorphic extension of an arbitrary solution of
 \begin{equation}\label{geswa}
 \{g,s\} = \frac{a}{2} - \frac{\esiz \alfa',\nu'\esde}{2} + i \, {\rm
 det}(\nu,\alfa,\nu',\alfa') + \{G,s\}.
 \end{equation}
Here, $\nu(s),G(s),\alfa(s)$ are defined following the previous geometric
construction.
\end{teo}
\begin{proof}
Given $\nu(s),\alfa(s):I\flecha \N^3$ as above, the maps $\beta(s)=
\frac{1}{2}\nu(s) +\alfa(s)$ and $V(s)= \frac{1}{2} \nu(s) -\alfa(s)$
describe a pair of Björling data, and thus generate a Bryant surface via
Theorem \ref{solbjbr}. Let $\phi$ be the solution to Liouville equation
for $c=1$ associate to such Bryant surface. Then $\phi(s,0)= -2 \esiz
\nu'(s),\nu'(s)\esde = a(s)$ and, by \eqref{datbjli}, $$\phi_t(s,0)= \esiz
\nu',\nu'\esde \left( \displaystyle\frac{(V\ex \nu')_0 +(V\ex
\nu')_3}{\nu_0 +\nu_3} + {\rm Im} \left(\frac{G''}{G'}\right) \right).$$
As $V\ex \nu' = \alfa \times \nu \times \nu' = \cF$, we deduce from
$\cF_0 +\cF_3 = \cL$ that $\phi_t (s,0)=d(s)$. In other words, $\phi$ is
the solution to \eqref{caulior1}.

On the other hand, we know from Umehara-Yamada's relation \eqref{uya} that $\phi$ is
expressed as \eqref{resceuno}, where $g(z)$ is an arbitrary solution of $\{g,z\} =
-2 q(z) + \{G,z\} $. But we already computed in the proof of Theorem \ref{solbjbr}
the quantities $q(z),G(z)$ in terms of the Björling data $\beta,V$, and these are
given by \eqref{bjh6}, \eqref{himabj}. Putting together these equations we obtain
that the differential equation in the real line describing $g$ in terms of
$\nu(s),\alfa(s)$ is given by \eqref{geswa}. This finishes the proof.
\end{proof}

The interaction of this result and Theorem \ref{curvas} yields the
following consequence, which reduces the problem of integrating the Frenet
equations in $\S^2$ to the one of finding a function with prescribed
Schwarzian derivative.
 \begin{cor}
Let $\gamma(s)$ be the only (up to congruences) curve in $\S^2$
parametrized by arclength, with geodesic curvature $\kappa(s)$. Then
$\gamma(s)= \pi^{-1} (g(s))$, where $\pi :\S^2\flecha \C\cup \{\8\}$ is
the stereographic projection and $g(s)$ is an arbitrary solution in
$\C\cup \{\8\} $ of \eqref{geswa}. Here, $\nu(s),G(s),\alfa(s)$ are defined following the previous geometric
construction for the choices $a(s)\equiv 1$, $d(s)= -2\kappa(s)$.
 \end{cor}

Finally, we make the following closing remark. We have used both geometric
and analytic methods to study the Cauchy problem for Liouville equation.
Nevertheless, the analytic approach is only applied to the
\emph{complexified} form of Liouville equation, and does not recover the
developing map of the solution. In contrast, when considering the
geometric method these two limitations disappear.

A phenomenon of this type was already noted by Liouville \cite{Lio}.
Indeed, rather that considering the elliptic case $\lap (\log \phi)= {\rm
const}\, \phi$, Liouville considered instead the hyperbolic case $(\log
\phi)_{uv} = {\rm const} \, \phi$, and provided an analytic resolution for
it. However, it comes clear from his comments that before giving this
analytic proof, he had already achieved a geometric resolution of the
elliptic equation, and that this solution was in essence the one obtained
in \cite{Lio} for the hyperbolic case by purely analytic methods.

In other words, when Liouville tried to provide an analytic resolution of
the elliptic equation $\lap (\log \phi)= {\rm const} \, \phi$, he found
the need to complexify this equation (which, once in this form, might be
seen as a hyperbolic real equation after the change $u=z$, $v= \bar{z}$).
In particular, with this technique the developing map $g$ is not
recovered.

Therefore, it seems that the geometric nature of Liouville equation gives
geometric methods a certain advantage with respect to the purely analytic
ones in its study.

\def\refname{References}

\hspace{0.4cm}

\noindent J.A. Gálvez was partially supported by MCYT-FEDER, Grant no. BFM2001-3318.

\hspace{0.2cm}

\noindent P. Mira was partially supported by MCYT, Grant no. BFM2001-2871 and CARM
Programa Séneca, Grant no PI-3/00854/FS/01.


\begin{thebibliography}{9}

\bibitem[ACM]{ACM} L.J. Al\' \i as, R.M.B. Chaves, P. Mira,  Bj\"orling problem for
maximal surfaces in Lorentz-Minkowski space, {\it Math. Proc. Cambridge
Philos. Soc.} {\bf 134} (2003), 289--316.


\bibitem[BPS]{Bob} A.I. Bobenko, T.V. Pavlyukevich, B.A. Springborn,
Hyperbolic constant mean curvature one surfaces: Spinor representation and
trinoids in hypergeometric functions, {\it Math. Z.} {\bf 245} (2003),
63--91.

\bibitem[Bry]{Bry} R.L. Bryant, Surfaces of mean curvature one in
hyperbolic space, {\it Astérisque}, {\bf 154-155} (1987), 321--347.

\bibitem[CHR]{CHR} P. Collin, L. Hauswirth, H. Rosenberg, The geometry of
finite topology Bryant surfaces, {\it Ann. of Math.} {\bf 153} (2001),
623--659.

\bibitem[DHKW]{DHKW}  U. Dierkes, S. Hildebrant, A. Küster, O. Wohlrab,
{\it Minimal Surfaces I}. Springer-Verlag, A series of comprehensive
studies in mathematics, 295, 1992.

\bibitem[EaTo]{EaTo} R. Sa Earp, E. Toubiana, On the geometry of constant mean curvature one
surfaces in hyperbolic space, {\it Illinois J. Math} {\bf 45} (2001),
371-401

\bibitem[GMM1]{GMM1} J.A. Gálvez, A. Martínez, F. Milán, Flat surfaces in the
hyperbolic $3$-space, {\it Math. Ann.}, {\bf 316} (2000), 419--435.

\bibitem[GMM2]{GMM2} J.A. Gálvez, A. Martínez, F. Milán, {\it Complete linear
Weingarten surfaces of Bryant type. A Plateau problem at infinity},
preprint.

\bibitem[GaMi1]{GaMi1} J.A. Gálvez, P. Mira, {\it Isometric
immersions of $\R^2$ into $\R^4$ and perturbation of Hopf tori}, preprint
(available at http://arXiv.org/math.DG/0301075)

\bibitem[GaMi2]{GaMi2} J.A. Gálvez, P. Mira, {\it Dense solutions to
the Cauchy problem for minimal surfaces}, preprint.

\bibitem[GaMi3]{GaMi3} J.A. Gálvez, P. Mira, {\it A Björling-type
problem for flat surfaces in hyperbolic $3$-space and its applications},
in preparation.

\bibitem[HRR]{HRR} L. Hauswirth, P. Roitman, H. Rosenberg, The
geometry of finite topology Bryant surfaces quasi-embedded in a hyperbolic
manifold, {\it J. Differential Geom.} {\bf 60} (2002), 55--101.

\bibitem[Kar]{Kar} H. Karcher, \emph{Hyperbolic constant mean curvature one
surfaces with compact fundamental domains}, preprint.

\bibitem[KUY]{KUY} M. Kokubu, M. Umehara, K. Yamada, {\it Flat fronts in hyperbolic $3$-space},
preprint.

\bibitem[Lio]{Lio} J. Liouville, Sur l'equation aux differences partielles
$\frac{\parc^2}{\parc u\parc v} \pm \frac{\landa}{2a^2} =0$,
{\it J. Math. Pures Appl.} {\bf 36} (1853), 71--72.

\bibitem[MiPa]{MiPa} P. Mira, J.A. Pastor, Helicoidal maximal surfaces
in Lorentz-Minkowski space, {\it Monatsh. Math.} (to appear).

\bibitem[Nit]{Nit} J.C.C. Nitsche, {\it Lectures on Minimal Surfaces},
Vol. I. Cambridge University Press, Cambridge, 1989.

\bibitem[Oss]{Oss} R. Osserman, {\it A survey of Minimal Surfaces}.
Enlarged republication of the 1969 original. Dover Publications, New York,
1986.

\bibitem[Sma]{Sma} A. Small, Surfaces of constant mean curvature $1$ in
$\H^3$ and algebraic curves on a quadric, {\it Proc. Amer. Math. Soc.}
{\bf 122} (1994), 1211--1220.

\bibitem[RUY1]{RUY1} W. Rossman, M. Umehara, K. Yamada, Irreducible
constant mean curvature $1$ surfaces in hyperbolic space with positive
genus, {\it Tohoku Math. J.} {\bf 49} (1997), 449--484.

\bibitem[RUY2]{RUY2} W. Rossman, M. Umehara, K. Yamada, \emph{Mean curvature $1$
surfaces with low total curvature in hyperbolic $3$-space I}, preprint.

\bibitem[RUY3]{RUY3} W. Rossman, M. Umehara, K. Yamada, Mean curvature $1$
surfaces with low total curvature in hyperbolic $3$-space II, {\it Tohoku
Math. J.} {\bf 55} (2003), 375--395

\bibitem[Ten]{Ten} K. Tenenblat, {\it Transformation of manifolds and applications
to differential equations,} Pitman Monographs and Surveys in Pure and
Applied Mathematics. Longman, Harlow, 1998.

\bibitem[UmYa1]{UY1} M. Umehara, K. Yamada, Complete surfaces of constant mean
curvature-$1$ in the hyperbolic $3$-space, {\it Ann. of Math.} {\bf 137}
(1993), 611--638.

\bibitem[UmYa2]{UY2} M. Umehara, K. Yamada, A parametrization of the
Weierstrass formulae and perturbation of complete minimal surfaces in
$\R^3$ into the hyperbolic $3$-space, {\it J. Reine Angew. Math.}
{\bf 432} (1992), 93--116.

\bibitem[UmYa3]{UYMath} M. Umehara, K. Yamada, Surfaces of constant mean
curvature-$c$ in $H^3(-c^2)$ with prescribed hyperbolic Gauss map, {\it
Math. Ann.} {\bf 304} (1996), 203--224.


\bibitem[Yu]{Yu} Z. Yu, The value distribution of hyperbolic Gauss maps,
{\it Proc. Amer. Math. Soc.} {\bf 125} (1997), 2997--3001.
\end{thebibliography}
\end{document}